\documentclass[11pt]{amsart}

\usepackage{graphicx, color}
\usepackage{amsmath, amssymb}
\usepackage{mathrsfs}
\usepackage{amsfonts}
\usepackage[all]{xy}
\usepackage{bbm}
\usepackage[margin=3cm]{geometry}

\numberwithin{equation}{section}

\DeclareMathOperator{\GL}{GL}
\DeclareMathOperator{\Aut}{Aut}
\DeclareMathOperator{\End}{End}

\DeclareMathOperator{\Ker}{Ker}
\DeclareMathOperator{\Diff}{Diff}

\DeclareMathOperator{\KZ}{KZ}

\def\lam{\lambda}
\def\Lam{\Lambda}
\def\endproof{\hfill{$\square$} }
\def\C{\mathbb{C}}
\def\R{\mathcal{R}}
\def\bR{\mathbb{R}}
\def\Q{\mathcal{Q}}
\def\Z{\mathbb{Z}}
\def\alp{\alpha}
\def\H{\mathcal{H}}
\def\iso{\xrightarrow{\sim}}
\def\L{\mathcal{L}}

\def\S{\mathbb{S}}
\def\bD{\mathbb{D}}
\def\rd{\partial}
\def\A{\mathbb{A}}
\def\M{\mathfrak{M}}
\def\CA{\mathcal{A}}
\def\TC{\widetilde{\mathcal{C}}}
\def\CC{\mathcal{C}}
\def\TA{\widetilde{\mathcal{A}}}
\def\TB{\widetilde{B}}
\def\g{\mathfrak{g}}
\def\one{\mathbbm{1}}

\def\Sec{\mathcal{S}}
\def\O{\mathcal{O}}

\def\d{\mathbf{d}}
\def\k{\mathbf{k}}
\def\N{\mathcal{N}}
\def\l{\mathbf{l}}
\def\bO{\mathbb{O}}

\newtheorem{thm}{Theorem}[section]
\newtheorem{prop}[thm]{Proposition}
\newtheorem{cor}[thm]{Corollary}
\newtheorem{lem}[thm]{Lemma}
\newtheorem{conj}[thm]{Conjecture}

\theoremstyle{definition}
\newtheorem{defi}[thm]{Definition}
\newtheorem{rem}[thm]{Remark}

\newcommand{\abv}[1]{\lvert #1 \rvert}
\newcommand{\br}[1]{\langle #1 \rangle}

\author{Akishi Ikeda}
\address{Kavli Institute for the Physics and Mathematics of the Universe (WPI), UTIAS, 
The University of Tokyo, Kashiwa, Chiba 277-8583, Japan}
\email{akishi.ikeda@ipmu.jp}


\begin{document}

\title[Representations of framed braid groups]
{Homological and monodromy representations of framed braid groups}

\begin{abstract}
In this paper, we introduce two new classes of representations of 
the framed braid groups. One is the homological representation 
constructed as the action of a mapping class group 
on a certain homology group. 
The other is 
the monodromy representation of the confluent KZ equation, which is 
a generalization of the KZ equation to have irregular singularities.   
We also give a conjectural equivalence between these two classes of 
representations. 
\end{abstract}

\maketitle

\section{Introduction}
The framed braid group $FB_n$ \cite{KS} is the semi-direct product 
$FB_n \cong \Z^n \rtimes B_n$ where the braid group 
$B_n$ acts on $\Z^n$ as permutations of 
components through the projection 
$B_n \to \mathfrak{S}_n$ on the symmetric group of degree $n$. 
Similar to the diagrammatic description of the braid group $B_n$ 
by using $n$ strings, the framed braid group $FB_n$ can be described 
graphically by using $n$ ribbons. The component $B_n$ describes crossings of ribbons 
and the component $\Z^n$ describes the number of twists of each ribbon. 
In particular, the plat closure of an 
element of  $FB_n$ gives a framed link. 
Whereas there are many studies of representations of the braid groups, 
representations of the framed braid groups are less known. 

The purpose of this paper is to introduce two new classes of representations of 
the framed braid groups. (These are also new classes as representations 
of the braid groups.)  
One is the homological representation 
constructed as the action of a mapping class on a certain homology group. 
This construction generalizes the homological representation of 
the braid group introduced by Lawrence \cite{Law}. 
The original Lawrence representation appears as the quotient 
of our representation. As special cases, our representations contain 
natural extensions of the reduced Burau representation \cite{Bur} and 
the Lawrence-Krammer-Bigelow representation \cite{Big,Kra1,Kra2} 
to the framed braid groups. 

The other is the monodromy representation of the confluent 
KZ equation 
\cite{JNS}, which is 
a generalization of the KZ equation \cite{KZ} to have irregular singularities. 
The confluent KZ equation also appears in 
two dimensional conformal field theory 
as the differential equation for irregular conformal 
blocks of the Wess-Zumino-Witten (WZW) model \cite{GL}. 
Thus our representation 
also describes the monodromy of irregular conformal blocks.

In \cite{Koh1}, Kohno showed that 
the Lawrence homological representation of the braid group 
is equivalent to the monodromy representation of the KZ equation on 
the space of singular vectors. 
Following his result, we also give a conjectural equivalence between 
the homological representation and the monodromy 
representation of the framed braid group. 

\subsection{Homological representations}
We summarize the construction of the homological representations of 
the framed braid groups from Section \ref{sec:homological}.
For a positive integer $n>0$, let $FB_n$ be the framed group 
of $n$ ribbons (see Definition \ref{def:framed_braid}).  
Set $\R=\Z[q^{\pm 1},t^{\pm 1}]$.
The homological representation
\[
\rho_{n,m}^{(r)} \colon FB_n \to \Aut_{\R} \H_{n,m}^{(r)}.
\]
is parametrized by positive integers $m,r>0$ where $\H_{n,m}^{(r)}$ 
is a certain relative homology group constructed as follows.
Let $\bD$ be a closed disk and take disjoint $n$ closed disks 
$D_1,\dots,D_n \subset \mathrm{Int} \bD$ from the interior of $\bD$. 
An open interval $A \subset \rd D_k$ is called a marked arc. 
Take disjoint  $r$ marked arcs from each boundary $\rd D_k$ and 
denote by $\A^{(r)}$ the set of such $rn$ marked arcs (Figure \ref{Fig:Arcs}). 
Define the surface $\S_n^{(r)}$ (Figure \ref{Fig:Surface}) by
\[
\S_n^{(r)}=\bD \setminus (D_1 \cup \cdots \cup D_n) \cup
\bigcup_{A \in \A^{(r)}}A.
\]
Let $\CC_{n,m}^{(r)}$ be the configuration space of unordered $m$ distinct 
points in $\S_n^{(r)}$:
\[
\CC_{n,m}^{(r)}=\{(t_1,\dots,t_m) \in 
(\S_n^{(r)})^m \,\vert \,t_i \neq t_j \text{ if }i \neq j \} 
\slash \mathfrak{S}_n.
\]
Then there is a group homomorphism from the fundamental group of $\CC^{(r)}_{n,m}$ 
to a free abelian group of rank two
\[
\alp \colon \pi_1(\CC^{(r)}_{n,m},\d) \to \br{q} \oplus \br{t}
\]
where the generator $q$ corresponds to 
the loop around the cylinders $\{t_1 \in D_k\}$ for $k=1,\dots,n$ 
and the generator $t$ corresponds to the loop around 
the hyperplanes $\{t_i=t_j\}$ for $1 \le i<j \le m$. 
Let $\pi \colon \TC_{n,m}^{(r)} \to \CC_{n,m}^{(r)}$ the covering space 
corresponding to $\alp$. 
Introduce the subset $\CA^{(r)} \subset \CC_{n,m}^{(r)}$ by
\[
\CA^{(r)}:=\{\,\{t_1,\dots,t_m\} \in \CC_{n,m}^{(r)}\,\vert\, t_1 \in A
 \text{ for some } A \in \A^{(r)}   \,\}
\]
and its inverse image $\TA^{(r)}=\pi^{-1} (\CA^{(r)}) \subset \TC_{n,m}^{(r)} $. 
The homological representation is constructed on the relative homology group
\[
\H_{n,m}^{(r)}=H_m(\TC_{n,m}^{(r)},\TA^{(r)};\Z ).
\]
We note that $\H_{n,m}^{(r)}$ has an $\R$-module 
structure coming from the action of the deck transformation group 
$\br{q} \oplus \br{t}$.
The mapping class group $\M(\S_n^{(r)})$ is defined to be 
the group of isotopy classes of orientation preserving 
diffeomorphisms on $\S_n^{(r)}$ which fix the boundary 
$\rd \bD$ pointwise. 
Then $\M(\S_n^{(r)})$ naturally acts on $\H_{n,m}^{(r)}$.  
In addition, there is an isomorphism $FB_n \cong \M(\S_n^{(r)})$ 
by Proposition \ref{prop:MCG_iso}. Our main result 
about the homological representation is the following. 
\begin{thm}
For positive integers $m,r>0$, there is a representation of 
the framed braid group
\[
\rho_{n,m}^{(r)} \colon FB_n \to \Aut_{\R} \H_{n,m}^{(r)}
\]
which is constructed as the action of 
the mapping class group $\M(\S_n^{(r)})$ on 
the relative homology group $\H_{n,m}^{(r)}$. 
This representation has the following properties:
\begin{itemize}
\item[(1)]
There is a subrepresentation 
$\L_{n,m}^{(r)} \subset \H_{n,m}^{(r)}$ which is a free $\R$-module 
of rank
\[
\binom{rn+n+m-2}{m}
\]
spanned by certain homology classes, 
called the standard multifork classes.
\item[(2)]
We have an equality 
\[
\L_{n,m}^{(r)} \otimes_{\R}\Q = \H_{n,m}^{(r)} \otimes_{\R}\Q.
\]
over the field $\Q=\mathbb{Q}(q,t)$. In particular, 
the standard multifork classes form a basis of $\H_{n,m}^{(r)}$ 
over $\Q$. 
\end{itemize}
\end{thm}
Details are given in Section \ref{sec:homological}. 
The original Lawrence representation \cite{Law} appears as follow. 
In $\L_{n,m}^{(r)}$, there is a natural subrepresentation  
$\N_{n,m}^{(r)} \subset \L_{n,m}^{(r)}$ of rank 
\[
\binom{rn+n+m-2}{m}-\binom{n+m-2}{m},
\]
and the quotient representation $\L_{n,m}=\L_{n,m}^{(r)} \slash \N_{n,m}^{(r)}$ 
is equivalent to the Lawrence representation \cite{Law} (also see  
a nice review of the Lawrence representation \cite[Section 3.1]{Ito}). 
In this construction, the factor $\Z^n$ of $FB_n \cong \Z^n \rtimes B_n$ acts on 
$\L_{n,m}$ trivially, and hence $\L_{n,m}$ has only information about 
the representation of the braid group $B_n$. 
As special cases, $\L_{n,1}$ is the reduced Burau representation \cite{Bur}
and $\L_{n,2}$ is the 
Lawrence-Krammer-Bigelow representation \cite{Big,Kra1,Kra2}.
The faithfulness of the representation $\L_{n,m}$ of $B_n$ for $m \ge 2$
by \cite{Big,Kra2,Zhe2} implies that the representation $\L_{n,m}^{(r)}$ of 
$FB_n$ is also faithful for $m \ge 2$ (see Corollary \ref{cor:faithful}).  
The polynomial invariants of framed links associated with 
the reduced Burau representations 
will be discussed in \cite{Ike} (see Remark \ref{rem:Alexander}).

\subsection{Monodromy representations}
We see the construction of the monodromy representations of 
the framed braid groups. Before going to the confluent KZ equation, 
we first recall some facts about 
the KZ equation and 
its monodromy representations. 
The Knizhnik-Zamolodchikov (KZ) equation was introduced in \cite{KZ} 
as the differential equation which is satisfied by 
correlation functions (conformal blocks) 
of the WZW model. Let $\g$ be a semi-simple Lie algebra and $M$ be a 
$\g$-module.  The KZ equation is an integrable 
differential equation for a $M^{\otimes n}$-valued 
function on the configuration space 
\[
X_n=\{(z_1,\dots,z_n) \in \C^n \vert z_i \neq z_j \text{ if }i \neq j\}
\]
with regular singularities along the divisors $\{z_i=z_j\}$.  
In addition, it is invariant under the action of $\mathfrak{S}_n$. 
Therefore 
we can associate a representation of 
the braid group $B_n \to \GL(M^{\otimes n})$ 
as the monodromy representation. 
The Kohno-Drinfeld theorem \cite{Dri,Koh2} describes 
the monodromy representation of the KZ equation as 
an $R$-matrix representation of the quantum group 
$U_q(\g)$. 
(The case $\g=\mathfrak{sl}_2$ with vector representations 
was studied in \cite{TK}.)  

Extensions of the KZ equation 
to have irregular singularities were studied in \cite{BK,FMTV,JNS}. 
In \cite{BK,FMTV}, the KZ equation with an irregular singularity of 
Poincar\'e rank $1$ at $\infty$ was introduced. 
In \cite{JNS}, they defined the confluent KZ equation which 
has irregular singularities of arbitrary 
Poincar\'e rank on each divisor $\{z_i=z_j\}$ and at $\infty$ 
in the case $\g=\mathfrak{sl}_2$. 
The confluent KZ equation is written by using Gaudin Hamiltonians 
with irregular singularities. The Gaudin model with irregular singularities 
was studied in \cite{FFT}. 
In two dimensional conformal field theory, the confluent KZ equation 
appears as the differential equation for 
irregular conformal blocks of the WZW model \cite{GL}. 

We briefly review the confluent KZ equation \cite{JNS}. 
Details are given in Section \ref{sec:monodromy}. 
In the rest of this section, we assume  $\g=\mathfrak{sl}_2$ with 
the standard basis $\{E,F,H\}$. For a positive integer $r>0$, 
introduce the truncated current Lie algebra 
$\g^{(r)}:=\g[t] \slash t^{r+1}\g[t]$
where $\g[t]=\g \otimes \C[t]$ is a Lie algebra with the bracket 
$[X \otimes t^m, Y \otimes t^n]:=[X,Y] \otimes t^{m+n}$. 
Write $X_p=X \otimes t^p$. 
Then 
we can define the confluent Verma module $M_{\lam}(\gamma)$ from the 
highest weight vector 
\[
E_i v_{\lam}(\gamma)=0,\quad H_0v_{\lam}(\gamma)=\lam 
v_{\lam}(\gamma),\quad H_i v_{\lam}(\gamma)=\gamma^{(i)}v_{\lam}(\gamma)\quad
(i=1,\dots,r)
\]
where $\lam \in \C$ 
is a weight of $H_0$ and 
$\gamma=(\gamma^{(1)},\dots,\gamma^{(r)}) \in \C^r$ with $\gamma^{(r)} \neq 0$ 
are weights of $H_1,\dots,H_r$. We call $\gamma^{(1)},\dots,\gamma^{(r)}$ 
movable weights. 
In the definition of the confluent KZ equation, 
movable weights $\gamma^{(1)},\dots,\gamma^{(r)}$ are regarded as 
variables of the equation.  
So we consider  the space of movable weights $B^{(r)}=\C^{r-1} \times \C^*$.
For $R=(r_1,\dots,r_n) \in \Z_{>0}^n$ and 
$\Lam=(\lam_1,\dots,\lam_n) \in \C^n$, consider the 
confluent Verma module bundle
\[
E^{(R)}_{\Lam}\to X_n \times B^{(r_1)}\times \cdots \times B^{(r_n)}
\]
whose fiber over a point $(z,\gamma_1,\dots,\gamma_n) \in X_n 
\times B^{(r_1)}\times \cdots \times B^{(r_n)}$ is the tensor product of the 
confluent Verma modules $M_{\lam_1}(\gamma_1)\otimes 
\cdots \otimes M_{\lam_n}(\gamma_n)$. 
Then there is an integrable connection $\nabla^{\KZ}$ on $E_{\Lam}^{(R)}$ 
depending on a complex parameter $\kappa \in \C^*$, 
called the confluent KZ connection. 
The confluent KZ connection $\nabla^{\KZ}$ has irregular singularities 
along the divisors $\{z_i=z_j\}$ and $\{\gamma_i^{(r_i)}=0\}$.
In the case $R=(r,\dots,r)=(r^n)$ and $\Lam=(\lam,\dots,\lam)=(\lam^n)$, 
the action of the symmetric group $\mathfrak{S}_n$ on 
both the total space $E_{\lam^n}^{(r^n)}$ and the base space 
$X_n \times (B^{(r)})^n$ is well-defined. In addition, the confluent KZ connection 
is $\mathfrak{S}_n$-invariant in this setting. 
Hence, the confluent KZ connection $\nabla^{\KZ}$ descends to 
an integrable connection on the quotient vector bundle
\[
E_{\lam^n}^{(r^n)} \slash \mathfrak{S}_n \to 
( X_n \times (B^{(r)})^n)\slash \mathfrak{S}_n.
\]
By Proposition \ref{prop:fundamental}, the fundamental group of 
the base space $( X_n \times (B^{(r)})^n)\slash \mathfrak{S}_n$ is 
isomorphic to the framed braid group $FB_n$. 
Hence we obtain a representation of $FB_n$ as the monodromy 
representation of the confluent KZ connection. 
Since the rank of $E_{\lam^n}^{(r^n)}$ is infinite, we consider the restriction of 
the confluent KZ connection on some finite rank subbundles.
For a positive integer $m$, there is a finite rank subbundle 
\[
S^{(r^n)}[n\lam-2m] \subset E_{\lam^n}^{(r^n)}
\]
consisting of singular vectors of weight $(n\lam-2m)$ 
such that the restriction of the confluent KZ connection $\nabla^{\KZ}$ 
on $S^{(r^n)}[n\lam-2m]$ is well-defined. 
Our main result 
about the monodromy representation is the following. 
\begin{thm}
For positive integers $m,r>0$ and complex parameters 
$(\lam,\kappa) \in \C \times \C^*$, there is a representation of 
the framed braid group
\[
\theta_{\lam,\kappa}^{(r)} \colon FB_n \to \Aut_{\C} S_{\Gamma}^{(r^n)}[n\lam-2m]
\]
which is constructed as the monodromy representation of 
the confluent KZ connection $\nabla^{KZ}$ on the vector bundle 
$S^{(r^n)}[n\lam-2m]$ where $S_{\Gamma}^{(r^n)}[n\lam-2m]$ 
is a complex vector space given by
the fiber over a point 
$(z,\Gamma) \in ( X_n \times (B^{(r)})^n)\slash \mathfrak{S}_n$. 
The dimension of the representation is given by
\[
\dim_{\C} S_{\Gamma}^{(r^n)}[n\lam-2m]=\binom{rn+n+m-2}{m}.
\]
\end{thm}
A new feature of this construction is that the confluent KZ connection 
has a non-trivial monodromy along the weight $\gamma^{(r)} \neq 0$ 
of $H_r \in \g^{(r)}$, and this monodromy action changes framings, i.e. it 
gives the action of the factor $\Z^n$ of $FB_n \cong \Z^n \rtimes B_n$.
In Section \ref{sec:monodromy}, we 
reformulate the confluent KZ equation \cite{JNS} in more geometric language  
and give a detailed study of it.

\subsection{Conjectural equivalences}
\label{conjecture}
In \cite{Koh1}, Kohno established the equivalence between 
the Lawrence homological representation of the braid group 
and the monodromy representation 
of the KZ equation. 
We expect the following generalization of his result to the framed braid group. 
\begin{conj}
There is an open dense subset $U \subset \C \times \C^*$ 
such that if $(\lam,\kappa) \in U$, then the monodromy representation 
$\theta_{\lam,\kappa}^{(r)}$ of $FB_n$ on the complex vector space 
$S_{\Gamma}^{(r^n)}[n\lam-2m]$ is equivalent to the homological representation 
$\rho_{n,m}^{(r)}$ of $FB_n$ on the 
complex vector space $\H_{n,m}^{(r)}\otimes_{\R}\C$ 
under the specialization 
\[
q=e^{\frac{2 \pi \sqrt{-1} \lam }{\kappa} },\quad 
t=-e^{ -\frac{2 \pi \sqrt{-1}}{\kappa}  }
\]
of variables of $\R=\Z[q^{\pm 1},t^{\pm 1}]$. 
\end{conj}
Note that our sign convention of exponent is 
opposite to the convention in \cite{Koh1}. 
This is due to the choice of the positive direction 
of $q$ and $t$ for our fork rules in Figure \ref{Fig:Rules}.
Kohno's proof in \cite{Koh1} is based on
\begin{itemize}
\item integral representations of solutions of the KZ equation by 
hypergeometic integrals over homology cycles 
with local system coefficients \cite{DJMM,SV},
\item linear independence of hypergeometric integral solutions 
by the determinant formula \cite{Var},
\item identification of the homological representation with 
the above homology group to construct 
hypergeometric integral solutions \cite{Koh1}.
\end{itemize}
To show Conjecture \ref{conjecture} similarly, 
we need to develop integral representations of solutions of the confluent KZ 
equation and the determinant formula. Integral representations of solutions
were studied in \cite{JNS,NS}. They represented the solution of the confluent 
KZ equation as an integral of variables $t_1,\dots,t_m$ of
the function $\Phi(t,z,\gamma)=\Phi_1(z,\gamma)\Phi_2(t,z,\gamma)$ 
where
\begin{align*}
\Phi_2(t,z,\gamma)=\prod_{1 \le a < b \le m}(t_a-t_b)^{2 \slash \kappa}
\prod_{a=1}^m \prod_{i=1}^n
\left\{(t_a-z_i)^{-\lam_i \slash \kappa} \exp 
\left(\frac{1}{\kappa}\sum_{p=1}^{r_i}\frac{1}{p}\frac{\gamma_i^{(p)}}{(t_a-z_i)^p}  
\right)\right\}
\end{align*}
and $\Phi_1(z,\gamma)$ is some function only depending on $z$ and $\gamma$ 
(for details, see \cite{JNS,NS}). Actually, the definition of the homological 
representation is motivated by the construction of appropriate 
homology cycles to integrate variables $t_1,\dots,t_m$ of 
the function $\Phi_2(t,z,\gamma)$. For convergence of the integration, 
we need to choose the direction in which $t_a$ approaches to $z_i$ carefully.  
Depending on the positive integer $r_i$, we can find $r_i$ sectors 
in $0<|t_a-z_i|<\epsilon$ where $\Phi_2(t,z,\gamma)$ decays 
exponentially as $t_a \to z_i$. These sectors essentially correspond to marked arcs 
in the construction of the homological representation. 
In two dimensional conformal field theory, the problem of the choice of 
this direction 
relates to a free field realization of 
irregular vertex operators (confluent primary fields) 
by using screening charges. For details, see \cite{GT,NS}. 
Our construction of the surface with marked arcs also motivated by 
the work \cite{HKK} which treats quadratic differentials with 
exponential singularities, and the name ``marked arc'' is taken 
from \cite{HKK}. Indeed in the case $m=1$, the function $\Phi_2$ can be 
identified with the (multi-valued) quadratic differential on $\C$ 
with exponential singularities at $z_1,\dots,z_n$. 
Finally also in the case $m=1$, the hypergeometric solution of 
the confluent KZ equation 
is equivalent to the hypergeometric function which was 
extensively studied by Haraoka in \cite{Hara}. In this case, he constructed 
appropriate integration cycles and also gave the determinant formula.

\subsection*{Acknowledgements}
This paper was written while the author was visiting Perimeter Institute 
for Theoretical Physics by JSPS Program for Advancing Strategic
International Networks to Accelerate the Circulation of Talented Researchers.
The author is grateful to host researcher Kevin Costello and 
stimulating research environment in  Perimeter Institute. 
This work is supported by World Premier International
Research Center Initiative (WPI initiative), MEXT, Japan and 
JSPS KAKENHI Grant Number JP16K17588.

\section{Framed braid groups}
\subsection{Definition of framed braid groups}
\begin{figure}[!ht]
\centering
\includegraphics[scale=1]{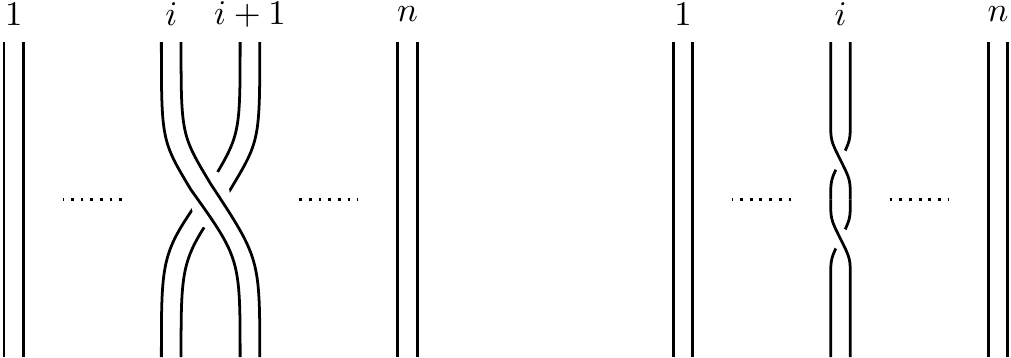}
\caption{Generators $\sigma_i$ and $\tau_i$.}
\label{Fig:Braid}
\end{figure}

Fix a positive integer $n>0$. 
Following \cite{KS}, we introduce the framed braid group as follows.
\begin{defi}
\label{def:framed_braid}
The  \textit{framed braid group} $FB_n$ is a group generated 
by $\sigma_1,\dots,\sigma_{n-1}$ and 
$\tau_1,\dots,\tau_n$ with the relations 
\begin{align*}
\sigma_i \sigma_{i+1}\sigma_i&=\sigma_{i+1}\sigma_i \sigma_{i+1}\\
\sigma_i \sigma_j &= \sigma_j \sigma_i \quad\text{if }|i-j| \ge 2 \\
\tau_i \tau_j&=\tau_j \tau_i\\
\sigma_i \tau_j&=
\begin{cases}
\,\tau_{i+1}\sigma_i &\text{if }j=i\\
\,\tau_i \sigma_i &\text{if }j=i+1\\
\,\tau_j \sigma_i &\text{if }j\neq i,i+1.
\end{cases}
\end{align*}
\end{defi}
Graphically generators of the framed braid group can be described
by using $n$ ribbons as in Figure \ref{Fig:Braid}. 
The generator $\sigma_i$ represents a crossing of the $i$-th and 
the $(i+1)$-st ribbons (left of Figure \ref{Fig:Braid}), 
and the generator $\tau_i$ represents a twisting of the $i$-th ribbon 
(right of Figure \ref{Fig:Braid}). 
The subgroup of $FB_n$ generated by $\sigma_1,\dots,\sigma_{n-1}$ 
becomes the \textit{braid group} $B_n$.
On the other hand, generators $\tau_1,\dots,\tau_n$ form a free 
abelian group of rank $n$ given by
\[
\left<\tau_1,\dots,\tau_n \right> \cong \Z^n, \quad
\tau_1^{l_1}\tau_2^{l_2} \cdots \tau_n^{l_n} \mapsto 
(l_1,l_2,\dots,l_n).
\]
Let $\mathfrak{S}_n$ be the symmetric group of degree $n$ and 
$\pi \colon B_n \to \mathfrak{S}_n$ be a surjective group homomorphism 
sending $\sigma_i$ to the transposition $(i,i+1)$. 
The kernel of $\pi$ is called the \textit{pure braid group} and denoted by $P_n$. 
The braid group $B_n$ acts 
on $\{1,\dots,n\}$ through the map $\pi$. We write this action as $\sigma(i)$ instead of 
$\pi(\sigma)(i)$. 
Then the framed braid group $FB_n$ is isomorphic to the semi-direct product 
$\Z^n \rtimes B_n$ where $B_n$ acts on $\Z^n$ by  
\[
\sigma \cdot (l_1,l_2,\dots,l_n):=(l_{\sigma(1)},l_{\sigma(2)},
\dots,l_{\sigma(n)}). 
\] 
The group homomorphism $\pi \colon B_n \to \mathfrak{S}_n$ can be extended to 
a group homomorphism $\pi^{\prime}\colon FB_n \to \mathfrak{S}_n$ 
by sending $\tau_i$ to the identity. 
We call the kernel of $\pi^{\prime}$ the \textit{pure framed braid group} 
and denote by $FP_n$. It is easy to see that $FP_n \cong \Z^n \times P_n$.

\subsection{Description as fundamental groups}
\label{sec:fundamental}
Let $X_n$ 
be the configuration space of ordered distinct $n$ points in $\C$, i.e.
\[
X_n:=\{\,(z_1,\dots,z_n) \in \C^n\,\vert\,z_i \neq z_j \text{ if }i \neq j\,\}.
\]
It is well-known that the fundamental group of 
$X_n$ is isomorphic to the pure braid group $P_n$. The quotient 
$X_n \slash \mathfrak{S}_n$ is the configuration space 
of unordered distinct $n$ points in $\C$ and its fundamental group 
is isomorphic to the braid group $B_n$. 
We give a similar description of the framed braid group. 
Let $T^n:=(S^1)^n$ be the $n$-dimensional torus. 
Then it is clear that the fundamental group of the direct product 
$T^n \times X_n$ is isomorphic to the pure framed braid group $FP_n \cong \Z^n \times P_n$.
Define the action of $s \in \mathfrak{S}_n$ on $T^n \times X_n$ by
\[
s \cdot (t_1,\dots,t_n,z_1,\dots,z_n):=(t_{s(1)},\dots,t_{s(n)},
z_{s(1)},\dots,z_{s(n)})
\]
for $(t_1,\dots,t_n) \in T^n$ and $(z_1,\dots,z_n) \in X_n$. 
In \cite{KS}, they showed the following. 
\begin{prop}[\cite{KS}, Proposition on page 2]
\label{prop:fundamental}
The quotient space
\[
(T^n \times X_n) \slash \mathfrak{S}_n
\] 
is a $K(FB_n, 1)$ space. Hence, its fundamental group is isomorphic to 
the framed braid group $FB_n$. 
\end{prop}
This result is used to construct the monodromy representations 
of the framed braid groups in Section \ref{sec:monodromy}. 

\subsection{Surfaces with marked arcs}
\label{sec:surface}
\begin{figure}[!ht]
\centering
\includegraphics[scale=0.9]{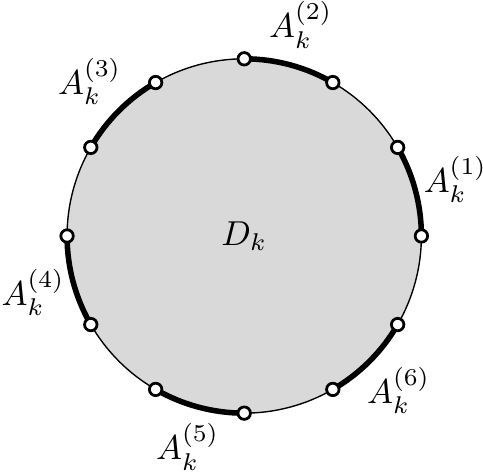}
\caption{Marked arcs in $\rd D_k$.}
\label{Fig:Arcs}
\end{figure}

Fix a positive integer $r$.  
Let $D_1,D_2,\dots,D_n$ be disjoint closed disks in $\C$ 
with centers  $1,2,\dots,n \in \C$ and the radius $1 \slash 4$, i.e. 
\[
D_k:=\left\{\,z \in \C\,\middle\vert\,  \abv{z-k} \le \frac{1}{4}   \,\right\}.
\]
Take a sufficiently large closed disk $\bD \subset \C$  
which contains $D_1,D_2,\dots,D_n$ in the interior. 
We define \textit{marked arcs}  
$A_k^{(s)} $ on $\rd D_k$  by 
\[
A_k^{(s)}:=\left\{\,k+\frac{1}{4}\exp\left(  \frac{i \pi \theta}{r} \right)
 \in \C\,\middle\vert\,\theta \in (2s-2,2s-1)  \,\right\}
\]
for $k=1,\dots,n$ and $s=1,\dots,r$ (see Figure \ref{Fig:Arcs}). 
Denote by $\A^{(r)}$ the set of these marked arcs:
\[
\A^{(r)}:=\left\{\,A_i^{(s)}\,
\middle\vert\,k=1,\dots,n \text{ and } s=1,\dots,r\,\right\}.
\]

\begin{figure}[!ht]
\centering
\includegraphics[scale=1]{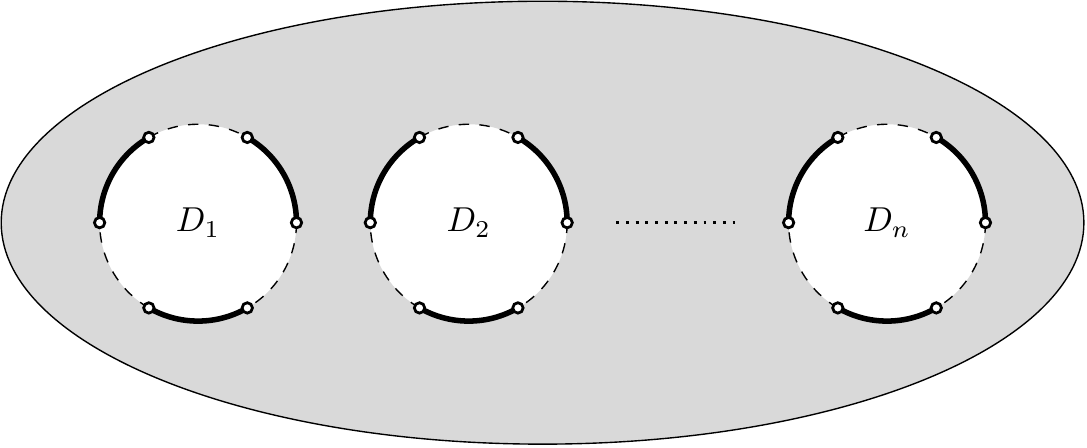}
\caption{Surface with marked arcs $\S_n^{(r)}$.}
\label{Fig:Surface}
\end{figure}

\begin{defi}
\label{defi:surface}
For positive integers $n>0$ and $r>0$, 
we define the surfaces $\S_n$ and $\S_n^{(r)}$ by
\[
\S_n:=\bD \setminus (D_1 \cup \cdots \cup D_n)
\]
and
\[
\S_n^{(r)}:=\S_n \cup \bigcup_{A \in \A^{(r)}}A.
\]
See Figure \ref{Fig:Surface}.
\end{defi}
The surface $\S_n^{(r)}$ is 
a connected oriented two dimensional smooth manifold with boundary 
\[
\rd \S_n^{(r)}=\rd \bD \cup \bigcup_{A \in \A^{(r)}}A,
\]
and non-compact. 
We call $\rd \bD$ the \textit{outer boundary}.

\subsection{Description as mapping class groups}
\label{sec:MCG}
Let $\S_n^{(r)}$ be the surface with marked arcs 
constructed in the previous section. 
Denote by $\Diff^+(\S_n^{(r)},\rd \bD)$ the group of orientation preserving 
diffeomorphisms on $\S_n^{(r)}$ which fix the outer boundary 
$\rd \bD$ pointwise. 
The subgroup of $\Diff^+(\S_n^{(r)},\rd \bD)$ consisting of 
diffeomorphisms which are isotopic to the identity is denoted 
by $\Diff^+_0(\S_n^{(r)},\rd \bD)$.

\begin{defi}
The \textit{mapping class group } $\M(\S_n^{(r)})$ is the quotient 
group defined by
\[
\M(\S_n^{(r)}):=\Diff^+(\S_n^{(r)},\rd \bD) \slash \Diff^+_0(\S_n^{(r)},\rd \bD).
\]
\end{defi}

\begin{figure}[!ht]
\centering
\includegraphics[scale=1]{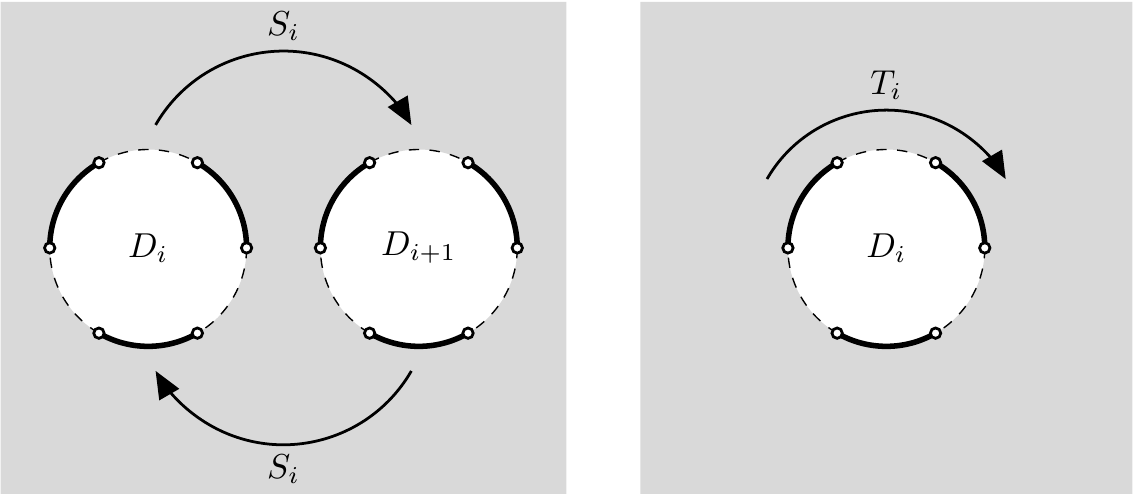}
\caption{Elements $S_i$ and $T_i$.}
\label{Fig:Twist}
\end{figure}

In the following, we see that the mapping class group $\M(\S_n^{(r)})$ 
is isomorphic to the framed braid group $FB_n$. 
We construct two classes of elements $S_1,\dots,S_{n-1} \in \M(\S_n^{(r)})$ 
and $T_1,\dots,T_n \in \M(\S_n^{(r)})$ as follows. 
\begin{itemize}
\item[(1)] The element $S_i \in \M(\S_n^{(r)})$ 
is a clockwise half twist  
between $\rd D_i$ and $\rd D_{i+1}$ which keeps the angle of 
circles $\rd D_i$ and $\rd D_{i+1}$, and hence  
interchanges arcs $A_i^{(s)}$ and $A_{i+1}^{(s)}$ 
(see the left of Figure \ref{Fig:Twist}).
\item[(2)] The element $T_i \in \M(\S_n^{(r)})$ 
is a clockwise rotation of the boundary $\rd D_i$ which 
maps $A_i^{(s)}$ to $A_i^{(s-1)}$ (see the right of Figure \ref{Fig:Twist}). 

\end{itemize}
\begin{prop}
\label{prop:MCG_iso}
There is an isomorphism of groups
\[
FB_n \iso \M(\S_n^{(r)})
\]
given by
\[
\sigma_i \mapsto S_i \quad
\text{and}\quad \tau_i \mapsto T_i.
\]
\end{prop}
\textbf{Proof. }Let $\M(\S_n)$ be the group of isotopy classes of 
orientation preserving diffeomorphisms on $\S_n$ which fix 
$\rd \bD$ pointwise. 
Note that $S_i$ is also well-defined on $\S_n \subset \S_n^{(r)}$. 
It is well-known that 
the correspondence
\[
B_n \to \M(\S_n),\quad \sigma_i \mapsto S_i
\]
gives an isomorphism between 
the braid group $B_n$ and the mapping class group 
$\M(\S_n)$. 
There is a natural surjection $\varphi \colon \M(\S_n^{(r)}) \to \M(\S_n)$ 
and it is easy to see that $\ker \varphi$ is generated by $T_1,\dots,T_n$. 
Hence we obtain a short exact sequence
\[
0 \to \left\langle T_1,\dots,T_n  \right\rangle \to \M(\S_n^{(r)}) 
\to \M(\S_n) \to 0
\]
and the group $\left\langle T_1,\dots,T_n  \right\rangle$ is isomorphic 
to $\Z^n$. 
\endproof
\medskip

This result is used to construct the homological representations 
of the framed braid groups in the next section.

\section{Homological representations}
\label{sec:homological}
\subsection{Configuration spaces and certain subsets}
Fix a positive integer $m>0$. Let $\S_n^{(r)}$ be 
the surface with marked arcs in Definition \ref{defi:surface}. 
Denote by $\CC_m(\S_n^{(r)})$ the configuration space 
of unordered distinct $m$ points in $\S_n^{(r)}$, i.e.
\[
\CC_m(\S_n^{(r)}):=\{\,(t_1,\dots,t_m) \in (\S_n^{(r)})^m\,\vert\,
t_i \neq t_j \text{ if }i\neq j  \,\} \slash \mathfrak{S}_m
\] 
where $\mathfrak{S}_m$ is the symmetric group of degree $m$. 
Write   
$\CC_{n,m}^{(r)}:=\CC_m(\S_n^{(r)})$ for simplicity. 
We use the notation $\{t_1,\dots,t_m\} \in \CC_{n,m}^{(r)}$ 
to represent unordered $m$ points. 
\begin{defi}  
The codimension one subset $\CA^{(r)} \subset \CC_{n,m}^{(r)}$  
is defined as a set of 
unordered distinct $m$ points in $\S_n^{(r)}$ 
such that at least one of $m$ points stays 
in some marked arc in $\A^{(r)}$, i.e.
\[
\CA^{(r)}:=\{\,\{t_1,\dots,t_m\} \in \CC_{n,m}^{(r)}\,\vert\, t_1 \in A
 \text{ for some } A\in \A^{(r)}   \,\}.
\]
\end{defi}
We introduce a specified base point in $\CC_{n,m}^{(r)}$ as follows. 
\begin{defi}
\label{def:basepoint}
Let $d_1,\dots,d_m$ be distinct $m$ points in the 
outer boundary $\rd \bD$ and 
assume that they lie in the lower half plane with the order 
$d_1,\dots,d_m$ from left to right (see Figure \ref{Fig:Fork}).   
We define the base point $\mathbf{d} \in \CC_{n,m}^{(r)}$ by 
$\d:=\{d_1,\dots,d_m\}$.
\end{defi}

\subsection{Relative homology groups}
In this section, we introduce certain relative homology groups on which 
our homological representations are constructed. 
Assume that $m \ge 2$. 
Then the first homology group of $\CC^{(r)}_{n,m}$ is given by
\[
H_1(\CC^{(r)}_{n,m};\Z) \cong \Z^{\oplus n} \oplus \Z
\]
where the first $n$ components correspond to the loops around the 
cylinders $\{|t_1-k| <(1 \slash 4)  \}$ for $k=1,\dots,n$, 
and the last component 
corresponds to the loop around the union of the 
hyperplanes $\{t_i=t_j\}$ for $1 \le i<j \le m$. 
Let $\br{q} \oplus \br{t}$ be a free abelian group of rank two 
with multiplicative generators $q$ and $t$.  
Define the group homomorphism 
\[
\alp \colon \pi_1(\CC^{(r)}_{n,m},\d) \to \br{q} \oplus \br{t} 
\]
by composing the map
\[
H_1(\CC^{(r)}_{n,m};\Z) \iso \Z^{\oplus n} \oplus \Z \to 
\br{q} \oplus \br{t} ,\quad (x_1,\dots,x_n,y) \mapsto (q^{x_1+\cdots+ x_n},t^y)
\]
and the abelianization map $\pi_1(\CC^{(r)}_{n,m},\d) 
\to H_1(\CC^{(r)}_{n,m};\Z)$. 
Let 
\[
\pi \colon \TC^{(r)}_{n,m}\to \CC^{(r)}_{n,m}
\]
the covering space corresponding to the normal subgroup 
$\Ker \alp \subset \pi_1(\CC^{(r)}_{n,m},\d) $.  
Define the subset
\[
\TA^{(r)}:=\pi^{-1}\left(\CA^{(r)} \right)  \subset \TC^{(r)}_{n,m} .
\]
The group $\br{q} \oplus \br{t}$ acts on $\TC_{n,m}^{(r)}$ and $\TA^{(r)}$ as 
deck transformations. Our homological representations are constructed 
on the following homology groups. Set $\R:=\Z[q^{\pm 1},t^{\pm 1}]$. 
\begin{defi}
An $\R$-module $\H_{n,m}^{(r)}$ is defined
to be the relative homology group
\[
\H_{n,m}^{(r)}:=H_{m}\left(\,\TC^{(r)}_{n,m}\,,\,\TA^{(r)};\Z \right)
\] 
where the $\R$-module structure is induced from the action of 
the deck transformation group $\br{q} \oplus \br{t}$.  
\end{defi}

\begin{rem}
In the case $m=1$, since $H_1(\CC_{n,1}^{(r)};\Z) \cong \Z^n$, 
the above constructions are considered by using 
the group homomorphism
$\alp \colon \pi_1(\CC^{(r)}_{n,1},\d) \to \br{q}$ instead. 
In the following sections, we always work over $\R=\Z[q^{\pm 1}]$ 
in the case $m=1$.
\end{rem}

\subsection{Construction of homological representations}
In this section, we define the action of the mapping class group $\M(\S_n^{(r)})$ 
on $\H_{n,m}^{(r)}$. Recall from Section \ref{sec:MCG} that 
the group $ \Diff^+(\S_n^{(r)},\rd \bD)$ consists of orientation preserving 
diffeomorphisms which fix $\rd \bD$ pointwise. 
A diffeomorphism $f \in \Diff^+(\S_n^{(r)},\rd \bD)$ 
induces a new diffeomorphism $\widehat{f} \colon \CC_{n,m}^{(r)} 
\to \CC_{n,m}^{(r)}$ 
defined by 
\[
\widehat{f}(\{t_1,\dots,t_m\}):=\{f(t_1),\dots,f(t_m)\}
\]
for $\{t_1,\dots,t_m\} \in \CC_{n,m}^{(r)}$. 
The base point $\mathbf{d}$ is fixed by $\widehat{f}$  
since $f$ fixes  $\rd \bD$ pointwise. 
There is a unique lift $\widetilde{f} \colon 
\TC_{n,m}^{(r)} \to \TC_{n,m}^{(r)}$ of  $\widehat{f}$ 
which fixes each point in $\pi^{-1}(\d)$ and commutes with 
the action of $\br{q} \oplus \br{t}$. In addition, $\widetilde{f}$ preserves 
the subset $\TA^{(r)} \subset \TC_{n,m}^{(r)}$ 
since $f$ preserves marked arcs $\A^{(r)}$. 
Hence $f$ induces an $\R$-linear map
\[
\widetilde{f}_* \colon \H_{n,m}^{(r)} \to \H_{n,m}^{(r)}
\]
on the relative homology group $\H_{n,m}^{(r)}$. 
It is clear that if $f$ is isotopic to the identity, 
then the induced linear map $\widetilde{f}_*$ is the identity. 
This implies that the action of the mapping class group $\M(\S_n^{(r)})$ 
on $\H_{n,m}^{(r)}$ is well-defined. Thus we obtain the following definition. 
\begin{defi}
Under the identification $FB_n \iso \M(\S_n^{(r)})$ in 
Proposition \ref{prop:MCG_iso}, 
the \textit{homological representation of the framed group}
\[
\rho_{n,m}^{(r)} \colon FB_n \to \Aut_{\R}\,\H_{n,m}^{(r)}
\]
is defined by $f \mapsto \widetilde{f}_*$
for $f \in FB_n$.
\end{defi}

\subsection{Multiforks}
\begin{figure}[!ht]
\centering
\includegraphics[scale=1.05]{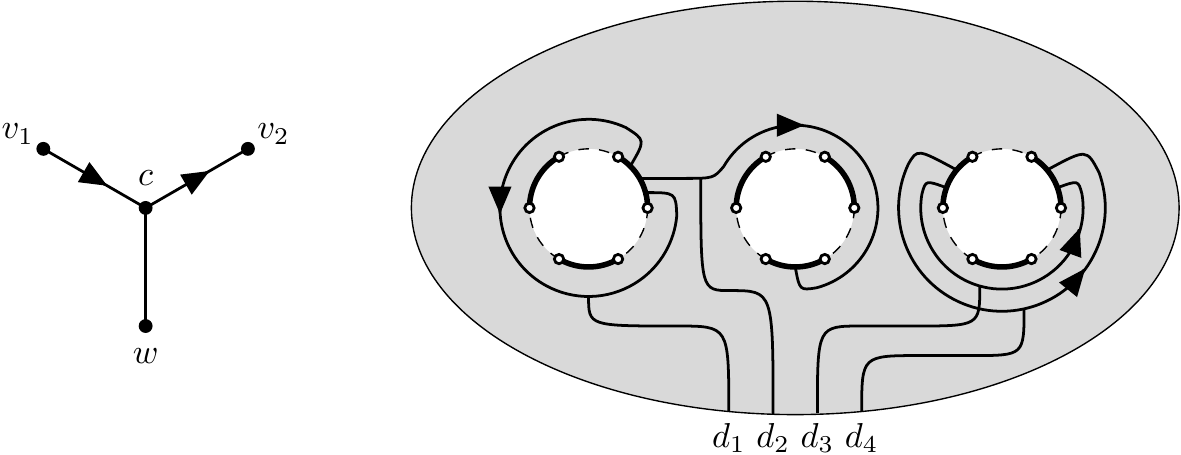}
\caption{The graph $Y$ and a multifork}
\label{Fig:Fork}
\end{figure}

\begin{figure}[!ht]
\centering
\includegraphics[scale=0.82]{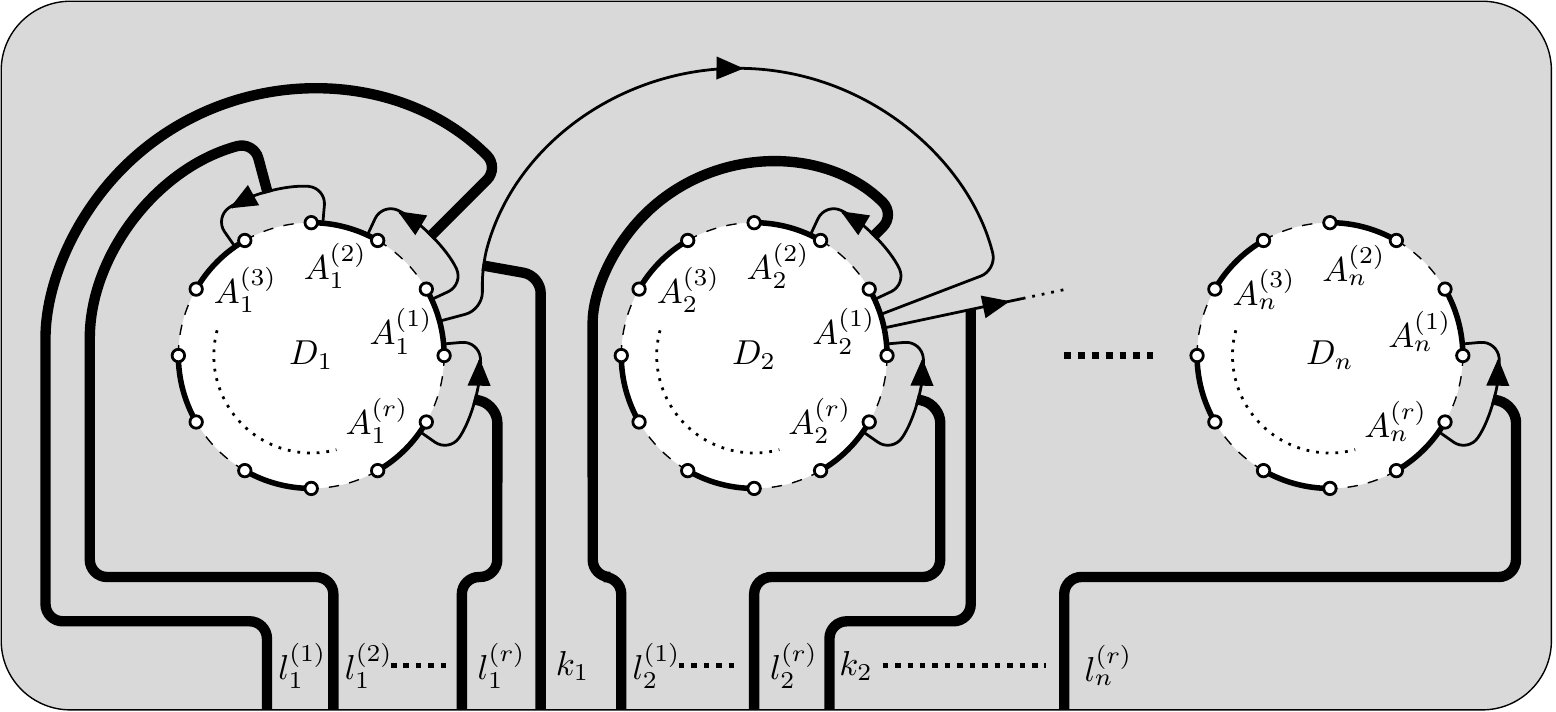}
\caption{Standard multifork}
\label{Fig:Standard}
\end{figure}

\begin{figure}[!ht]
\centering
\includegraphics[scale=1]{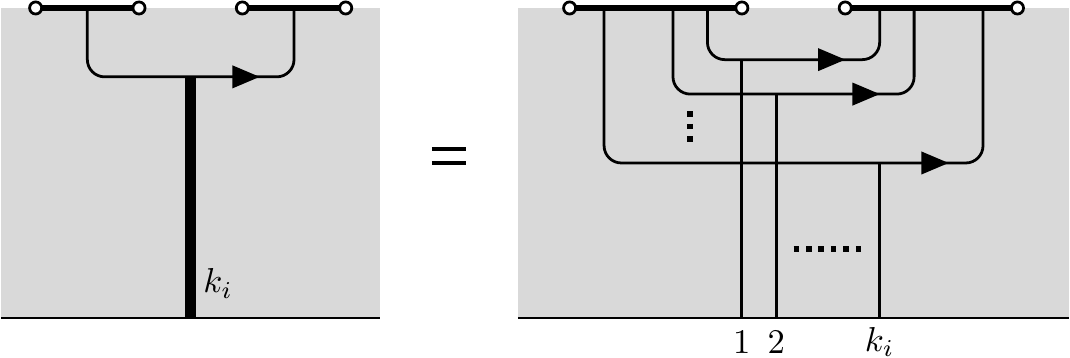}
\caption{Parallel forks}
\label{Fig:Parallel}
\end{figure}

\begin{figure}[!ht]
\centering
\includegraphics[scale=0.9]{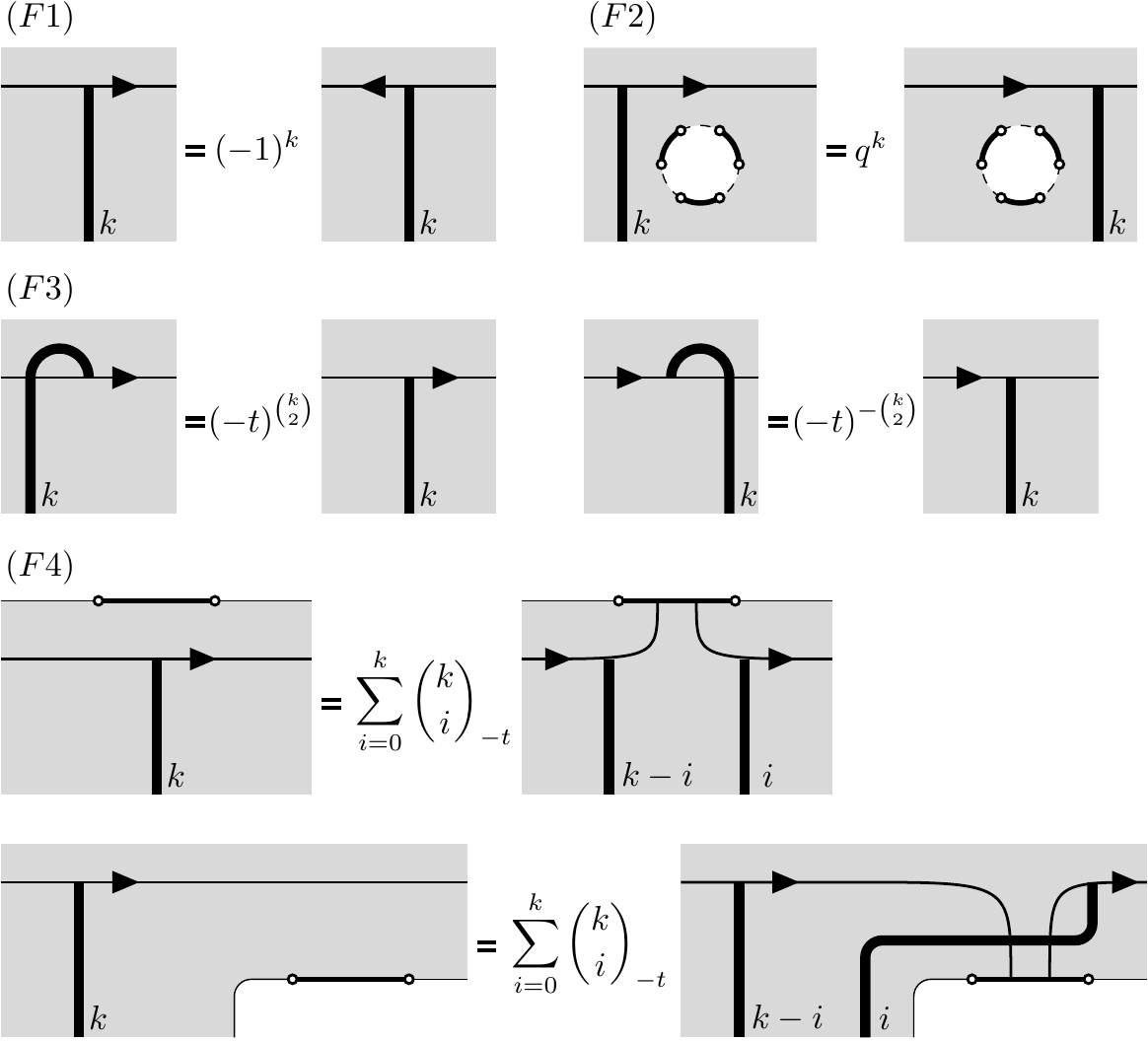}
\caption{Fork rules}
\label{Fig:Rules}
\end{figure}

We recall the notion of forks from \cite{Big,Kra1} and multiforks from \cite{Ito,Zhe}. 
We borrow the notation from \cite[Section 3]{Ito}. 
Let $\d=\{d_1,\dots,d_n\} \in \CC_{n,m}^{(r)}$ be the base point in 
Definition \ref{def:basepoint}. 
\begin{defi}
Let $Y$ be a graph consisting of 
four vertices $\{v_1,v_2,w, c\}$ and 
three edges $\{[v_1,c],[c,v_2],[w,c]  \}$ as in the left of Figure \ref{Fig:Fork}. 
We orient edges $[v_1,c]$ and $[c,v_2]$ as in Figure \ref{Fig:Fork}. 
A \textit{fork $F$ based on $d_i$} is the image of an embedding 
$\phi \colon Y \to \S_n^{(r)}$ satisfying the followings:
\begin{itemize}
\item $\phi(Y \setminus \{v_1,v_2,w\})$ lies in the interior of $\S_n^{(r)}$,
\item $\phi(w)=d_i$,
\item $\phi(v_1)  \in A_1 $ and $\phi(v_2)  \in A_2 $ 
for some $A_1,A_2 \in \A^{(r)} $, 
\item for $[v_1,v_2]:=[v_1,c] \cup [c,v_2]$, 
the closed arc $\phi([v_1,v_2])$ is not homotopic to 
a closed subinterval of  some $A \in \A^{(r)}$.
\end{itemize}
\end{defi}
The image $\phi([w,c])$ is called the \textit{handle} of $F$ and denoted by $H(F)$. 
The image $\phi([v_1,v_2])$ is called the \textit{tine edge} of $F$ and denoted by 
$T(F)$.
\begin{defi}
A \textit{multifork} is a family of forks $\mathbb{F}:=(F_1,\dots,F_m)$ 
satisfying the followings:
\begin{itemize}
\item $F_i$ is a fork based on $d_i$,
\item $T(F_i) \cap T(F_j) =\emptyset$ for $i \neq j$,
\item $H(F_i) \cap H(F_j) =\emptyset$ for $i \neq j$.
\end{itemize}
\end{defi} 
The right of Figure \ref{Fig:Fork} is an example of a multifork. 
For a multifork $\mathbb{F}=(F_1,\dots,F_m)$, 
let $\gamma_i \colon [0,1] \to \S_n$ be the path corresponding 
to the handle $H(F_i)$ with $\gamma_i(0)=d_i$.
Since all handles of $\mathbb{F}$ are disjoint, the path
\[
H(\mathbb{F}) \colon [0,1] \to \CC^{(r)}_{n,m}, \quad 
t \mapsto \{\gamma_1(t),\dots,\gamma_n(t)\}
\]
is well-defined. Note that $H(\mathbb{F})(0)=\d $. 
For the covering $\pi \colon \TC^{(r)}_{n,m}\to \CC^{(r)}_{n,m} $, 
fix a lift of the base point $\widetilde{\d} \in \pi^{-1}(\d)$. 
Then we can take a unique lift $
\widetilde{H(\mathbb{F})} \colon [0,1] \to \TC^{(r)}_{n,m}$ 
satisfying $\widetilde{H(\mathbb{F})}(0)=\widetilde{\d}$. 
For a given multifork $\mathbb{F}$, 
define the $m$ dimensional submanifold of $\CC^{(r)}_{n,m}$ by
\[
\Sigma(\mathbb{F}):=\{\,(t_1,\dots,t_m) \in\CC^{(r)}_{n,m} \,\vert\,\, 
t_i \in T(F_i)\,\}. 
\]
By definition, the boundary of $\Sigma(\mathbb{F})$ 
is contained in $\CA^{(r)}$. 
Let $\widetilde{\Sigma}(\mathbb{F})\subset \TC^{(r)}_{n,m}$ be the connected component 
of $\pi^{-1}(\Sigma(\mathbb{F}))$ containing 
the point $\widetilde{H(\mathbb{F})}(1)$. 
Since the boundary of $\widetilde{\Sigma}(\mathbb{F})$ 
is contained in $\TA^{(r)}$, it defines a homology class 
$[\widetilde{\Sigma}(\mathbb{F})] \in \H_{n,m}^{(r)}$, 
called the \textit{multifork class}.  We denote it by $[\mathbb{F}]$ 
instead of $[\widetilde{\Sigma}(\mathbb{F})]$. 
Introduce the set
\[
K_{n,m}^{(r)}:=\left\{\k=\left(k_1,\dots,k_{n-1},l_1^{(1)},
\dots,l_1^{(r)},\dots,l_n^{(1)},\dots,l_n^{(r)}\right)
\in  \Z_+^{n-1}\times (\Z_+^r)^n \middle\vert |\k|=m \right\}
\]
where 
\[
|\k|=\sum_{i=1}^{n-1}\,k_i+\sum_{i=1}^n \sum_{s=1}^{r} \,l_i^{(s)}.
\]
For $\k \in K_{n,m}^{(R)}$, we define the \textit{standard multifork} 
$\mathbb{F}_{\k}$ as in Figure \ref{Fig:Standard}. 
In Figure \ref{Fig:Standard}, a fork with a thick handle equipped with an integer 
$k_i$ or $l_i^{(s)}$ represents 
parallel $k_i$ or $l_i^{(s)}$ forks connecting $A_i^{(1)}$ and $A_{i+1}^{(1)}$, 
or $A_i^{(s)}$ and $A_i^{(s+1)}$ as in Figure \ref{Fig:Parallel}. 
By definition, the standard multifork classes $\{\,[\mathbb{F}_{\k}] 
\,\vert\, \k \in K_{n,m}^{(r)}\,\}$ are linearly independent in $\H_{n,m}^{(r)}$. 
In the relative homology group $\H_{n,m}^{(r)}$, multifork classes satisfy 
the diagrammatic formulas, called the \textit{fork rules} 
as in Figure \ref{Fig:Rules}. 
In Figure \ref{Fig:Rules}, we use the $(-t)$-binomial coefficient 
defined by
\[
\binom{k}{i}_{-t}:=\frac{[k]_{-t}\,!}{ [k-i]_{-t}\,!\,[i]_{-t}\,!  }
\]
where 
\[
[k]_{-t}\,!:=[1]_{-t}[2]_{-t} \cdots [k]_{-t}
\quad\text{and}\quad
[n]_{-t}:=\frac{(-t)^n-1}{(-t)-1}.
\]
By using these formulas, we can compute the action of $FB_n$ on 
multifork classes explicitly. In particular, any multifork class 
can be written as a linear combination of the standard multifork classes over $\R$. 
Hence we have the following. 
\begin{prop}
\label{prop:standard}
Let $\mathcal{L}_{n,m}^{(r)} \subset \H_{n,m}^{(r)}$ be a free $\R$-module of 
rank
\[
\left|K_{n,m}^{(r)} \right|=\binom{rn+n+m-2}{m}
\]
spanned by the standard multifork classes $\{\,[\mathbb{F}_{\k}] 
\,\vert\, \k \in K_{n,m}^{(r)}\,\}$. 
Then $\mathcal{L}_{n,m}^{(r)}$ is a representation of $FB_n$ over $\R$.  
\end{prop}
We also call $\L_{n,m}^{(r)}$ the homological representation. 
In Section \ref{sec:dimension}, we show that $\L_{n,m}^{(r)}$ coincides with 
$\H_{n,m}^{(r)}$ over the field $\mathbb{Q}(q,t)$, and therefore 
the standard multifork classes form 
a basis of $\H_{n,m}^{(r)}$ over $\mathbb{Q}(q,t)$. 

\subsection{Lawrence representations as quotients}
In this section, we see that the Lawrence representation \cite{Law} 
appears as the quotient of our homological representation. 
Define the set $K_{n,m}$ by
\[
K_{n,m}:=\left\{\k=\left(k_1,\dots,k_{n-1}\right)
\in \Z_+^{k-1}\,\middle\vert\, \sum_{i=1}^{k-1} k_i=m \,\right\}
\]
and regard $K_{n,m}$ as a subset of $K_{n,m}^{(r)}$. 
Let $\N_{n,m}^{(r)} \subset \L_{n,m}^{(r)}$ be the free $\R$-submodule of rank 
\[
\binom{rn+n+m-2}{m}-\binom{n+m-2}{m}
\]
spanned by the subset $\{\,[\mathbb{F}_{\k}] 
\,\vert\, \k \in K_{n,m}^{(r)} \setminus K_{n,m} \,\}$ of 
the standard multifork classes. 
Then we can check that $\N_{n,m}^{(r)}$ is closed under the 
action of $FB_n$ by using the fork rules in Figure \ref{Fig:Rules}. 
Hence the quotient representation $\L_{n,m}:=\L_{n,m}^{(r)} \slash \N_{n,m}^{(r)}$ 
is well-defined. 
\begin{prop}
The quotient $\L_{n,m}$ is a free $\R$-module of rank
\[
\binom{n+m-2}{m}
\]
spanned by $\{\,[\mathbb{F}_{\k}] \,\vert\, \k \in K_{n,m} \,\}$.  
The action of the factor $\Z^n \subset FB_n$ on 
the quotient representation $\L_{n,m}$ is trivial. 
As a representation of $B_n \subset FB_n$, the quotient 
representation $\L_{n,m}$ is equivalent to the Lawrence representation. 
\end{prop}
\textbf{Proof. }By the fork rules in Figure \ref{Fig:Rules}, 
it is easy to check that $\Z^n$ acts trivially on 
$\{\,[\mathbb{F}_{\k}] \,\vert\, \k \in K_{n,m} \,\}$ modulo $\N_{n,m}^{(r)}$. 
The basis $\{\,[\mathbb{F}_{\k}] \,\vert\, \k \in K_{n,m} \,\}$ and 
the action of $B_n$ on it modulo $\N_{n,m}^{(r)}$ 
precisely correspond to the description of 
the Lawrence representation in \cite{Ito,Zhe} by using multiforks. 
\endproof

\begin{cor}
\label{cor:faithful}
If $m \ge 2$, then the homological representation 
\[
\rho_{n,m}^{(r)} \colon FB_n \to \Aut_{\R}\,\L_{n,m}^{(r)}
\]
is faithful. 
\end{cor}
\textbf{Proof. }By definition, it is easy to see that 
the factor $\Z^n \subset FB_n$ 
acts faithfully on $\L_{n,m}^{(r)}$. On the other hand, the factor $B_n \subset FB_n$ 
acts faithfully on the quotient representation $\L_{n,m}$ 
by \cite{Big,Kra2} in the case $m=2$ and 
by \cite{Zhe2} in the case $m \ge 3$. 
\endproof

\subsection{Dimension of homological representations}
\label{sec:dimension}
Let $\Q:=\mathbb{Q}(q,t)$ be the quotient field of $\R=\Z[q^{\pm 1},t^{\pm 1}]$. 
In this section, we compute the dimension of the relative homology 
group $H_*\left(\,\TC^{(r)}_{n,m}\,,\,\TA^{(r)};\Z \right)$ over $\Q$. 
In particular, we recall the definition 
\[
\H_{m,n}^{(r)}=H_{m}\left(\,\TC^{(r)}_{n,m}\,,\,\TA^{(r)};\Z \right).
\]
\begin{thm}
\label{thm:dimension}
The dimension of the relative homology group 
 is given by
\[
\dim_{\Q}\,\H_{m,n}^{(r)}\otimes_{\R} \Q 
= \binom{rn+n+m-2}{m} 
\]
and
\[
\dim_{\Q}\,H_{k}\left(\,\TC^{(r)}_{n,m}\,,\,\TA^{(r)};\Z \right)\otimes_{\R} \Q = 0
\]
for $k \neq m$. 
\end{thm}
In order to show Theorem \ref{thm:dimension}, 
we begin by preparing some notations. 
For a topological space $X$, we denote by $\CC_{m}(X)$ the 
configuration space of unordered distinct $m$ points in $X$:
\[
\CC_{m}(X):=\{  \,(t_1,\dots,t_m) \in X^m\,\vert\,
t_i \neq t_j \text{ if }i \neq j   \,   \} \slash \mathfrak{S}_m.
\] 
Define a sequence of subsets 
\[
\emptyset=B_{-1}\subset 
B_{0}\subset B_{1}\subset 
\cdots \subset B_{m-p}\subset \cdots
 \subset B_{m-1} \subset B_{m}  \subset \CC_{n,m}^{(r)}
\]
by
\[
B_{m-p}\!:=\!
 \left\{\{t_1,\dots,t_m\} \in \CC_m (\S_n^{(r)} \setminus \rd \bD)  
 \,\middle\vert\, t_1 \in A_1,\dots,
t_p \in A_p \text{ for some } A_1,\dots,A_p \in \A^{(r)}   \right\}.
\]
In other words, at least $p$ points of 
$\{t_1,\dots,t_m\} \in \CC_m(\S_n^{(r)} \setminus \rd \bD)$ lie 
in some arcs of $\A^{(r)}$. 
\begin{rem}
\label{rem:corner}
The space 
\[
B_m=\bigsqcup_{p=0}^m (B_p \setminus B_{p-1})
\]
has a structure of a $2m$-dimensional smooth 
manifold with corners 
induced from the smooth structure on $\S_n^{(r)} \setminus \rd \bD$. 
Codimension $p$ corners of $B_m$
are given by $B_{m-p} \setminus B_{m-p-1}$, i.e. for any point 
$x \in B_{m-p} \setminus B_{m-p-1}$, there is an open subset $x \in U \subset B_m$ 
and a smooth map $f \colon U \to V$ on an open subset
$V \subset [0,\infty)^{p} \times \bR^{2m-p}$
such that $f$ is a diffeomorphism and $f(x)=0$. Since 
$[0,\infty)^{p} \times \bR^{2m-p}$ is homeomorphic to 
$[0,\infty) \times \bR^{2m-1}$ for $p \ge 1$, if we forget the smooth structure 
on $B_m$ and only consider the  
topological structure, then $B_m$ is a $2m$-dimensional 
topological manifold with boundary $\rd B_m=B_{m-1}$. In the following, 
our computation only depends on the topological structure on $B_m$. 
\end{rem}

We note that 
the pair $(B_{m},B_{m-1})$ is homotopy equivalent to 
the pair $(\CC_{n,m}^{(r)},\CA^{(r)})$ since 
$\S_n^{(r)} \setminus \rd \bD$ is homotopy equivalent to $\S_n^{(r)}$. 

For a subspace $U \subset \CC_{n,m}^{(r)}$, 
denote by $\widetilde{U}:=\pi^{-1}(U)$ its inverse image  via the covering map 
$\pi \colon \TC_{n,m}^{(r)} \to \CC_{n,m}^{(r)}$. 
Then the pair  $(\TB_m,\TB_{m-1})$ is also 
homotopy equivalent to 
the pair $(\TC_{n,m}^{(r)},\TA^{(r)})$. 
Therefore it is sufficient to compute 
the dimension of the relative homology group
\[
H_*(\,\TB_{m}\,,\,\TB_{m-1};\Z )\otimes_{\R} \Q.
\]
In the following of this section, 
we consider all homology groups over $\Q$ and 
write $H_*(\,-\,)$ instead of $H_*(\,-\,;\Z)\otimes_{\R}\Q$ for simplicity.

\begin{lem}
\label{lem:dim1}
The space $\TB_{p} \setminus \TB_{p-1}$ is 
an $(m+p)$-dimensional (smooth) manifold without boundary. 
The dimension of the homology group is given by
\[
\dim_{\Q} H_p (\TB_p \setminus \TB_{p-1})=
\binom{rn+m-p-1}{m-p} \binom{n+p-2}{p}  
\]
and 
\[
\dim_{\Q} H_k (\TB_p \setminus \TB_{p-1})=0
\]
for $k \neq p$.
\end{lem}
\textbf{Proof. }Set 
\[
\bO_n^{(r)}:=\S_n^{(r)} \setminus \S_n=\bigsqcup_{A \in \A^{(r)}}A.
\]
The space $\bO_n^{(r)}$ is an one dimensional manifold 
consisting of disjoint $rn$ open intervals. By definition, 
there is an isomorphism
\[
B_p \setminus B_{p-1} \cong \CC_{m-p}(\bO_n^{(r)}) 
\times \CC_p(\S_n \setminus \rd \bD),
\]
and hence $B_p \setminus B_{p-1}$ is an 
$(m+p)$-dimensional manifold without boundary. 
Since the covering map $\pi $ is a local homeomorphism, 
$\TB_p \setminus \TB_{p-1}=\pi^{-1}(B_p \setminus B_{p-1})$ 
is also an $(m+p)$-dimensional manifold without boundary. 
Let $L_{n,m-p}^{(r)}$ be a set
\[
L_{n,m-p}^{(r)}:=\left\{\,\l=\left(l_1^{(1)},
\dots,l_1^{(r)},\dots,l_n^{(1)},\dots,l_n^{(r)}\right)
\in   (\Z_+^r)^n\, \middle\vert\, \sum_{i=1}^n \sum_{s=1}^{r} \,l_i^{(s)}=m-p\,\right\}
\]
and for $\l \in L_{n,m-p}^{(r)}$, define a subset 
$O[\l]\subset \CC_{m-p}(\bO_n^{(r)})$ 
by
\[
O[\l] \!:=\!
 \left\{\{t_1,\dots,t_{m-p}\} \in \CC_{m-p}(\bO_n^{(r)})\,\middle\vert
\begin{array}{l}  
t_1,\dots,t_{l_1^{(1)}} \in A_1^{(1)},t_{l_1^{(1)}+1},\dots,t_{l_1^{(1)}+l_1^{(2)}}\in A_1^{(2)}, 
\dots,\! \\
 t_{l_1^{(1)}+\cdots+l_n^{(r-1)}+1},\dots,t_{l_1^{(1)}+\cdots+l_n^{(r)}} \in A_n^{(r)}
\end{array}
\,\right\}.
\]
Namely, the space $O[\l] $ consists of configurations 
in $\bO_n^{(r)}$ which satisfy that $l_i^{(s)}$ points of 
$\{t_1,\dots,t_{m-p}\} \in \CC_{m-p}(\bO_n^{(r)})$ lie in $A_i^{(s)}$. 
Then $\CC_{m-p}(\bO_n^{(r)})$ is decomposed into connected components
\[
\CC_{m-p}(\bO_n^{(r)})=\bigsqcup_{\l \in L_{n,m-p}^{(r)}}O[\l]
\]
and each $O[\l]$ is contractible by definition. 
On the other hand, it is known that 
\[
\dim_{\Q}H_p (\TC_p(\S_n \setminus \rd \bD))=\binom{n+p-2}{p}
\]
and $H_k(\TC_p(\S_n \setminus \rd \bD))=0$ for $k \neq p$ 
since it is the dimension of the Lawrence representation 
(see \cite[Section 2 and 3]{Koh1}). 
Since 
\[
\TB_p \setminus \TB_{p-1}\cong 
\bigsqcup_{\l \in L_{n,m-p}^{(r)}}(O[\l] \times \TC_p(\S_n \setminus \rd \bD))
\]
and 
\[
\left|L_{n,m-p}^{(r)} \right|=\binom{rn+m-p-1}{m-p},
\]
the result follows.
\endproof
\medskip

In the following computation, we use the cohomology 
with compact support, and Poincar\'e and Lefschetz dualities for 
non-compact manifolds. For these materials, 
we refer to \cite[Section 3.3]{Hat} 
and \cite[Section III]{Ive}. For a (locally compact) space $X$, 
we denote by $H^*_c(X;\Z)$ 
the singular cohomology with compact support. 
Since our spaces always admit the action of $\br{q} \oplus \br{t}$, 
we work over $\Q$ and write $H^*_c(\,-\,)$ instead of 
$H^*_c(\,-\,;\Z) \otimes_{\R}\Q$.

\begin{lem}
\label{lem:dim2}
We have the equations of dimensions
\[
\dim_{\Q} H^m_c(\TB_p)=\dim_{\Q} H^m_c(\TB_p \setminus \TB_{p-1})+
\dim_{\Q} H^m_c(\TB_{p-1})
\]
and 
\[
\dim_{\Q} H^k_c(\TB_p )=0
\]
for $k \neq m$.
\end{lem}
\textbf{Proof. }Since $\TB_p \setminus \TB_{p-1}$ is 
an $(m+p)$-dimensional manifold without boundary by Lemma \ref{lem:dim1},  
we have an isomorphism  
$H^k_c(\TB_p \setminus \TB_{p-1}) \cong H_{m+p-k}(\TB_p \setminus \TB_{p-1})$ 
by Poincar\'e duality (\cite[Theorem 3.35]{Hat}).  Again by Lemma \ref{lem:dim1}, 
$ H^k_c(\TB_p \setminus \TB_{p-1}) \neq 0$ if and only if $k=m$. 
Since $\TB_p \setminus \TB_{p-1}$ is an open subset of $\TB_p$ and 
$\TB_p - (\TB_p \setminus \TB_{p-1})=\TB_{p-1}$, 
we have a long exact sequence of cohomology with compact support 
(\cite[Section III.7]{Ive})
\[
\cdots \to  H^{k-1}_c(\TB_{p-1})\to H^k_c(\TB_p \setminus \TB_{p-1}) 
\to H^k_c(\TB_p) \to 
H^k_c(\TB_{p-1}) \to H^{k+1}_c(\TB_p \setminus \TB_{p-1})\to \cdots.
\]
We show the result by induction for $p$. First we consider the case $k \neq m$. 
By the long exact sequence, if 
$H^k_c(\TB_p \setminus \TB_{p-1})=0$ and $H^k_c(\TB_{p-1})=0$, 
then $H^k_c(\TB_p)=0$. By Lemma \ref{lem:dim1}, 
$H^k_c(\TB_0)=H^k_c(\TB_0 \setminus \TB_{-1})=0$ and 
$H^k_c(\TB_p \setminus \TB_{p-1})=0$ for $p\ge 0$ and $k \neq m$. 
Inductively, we have 
$H^k_c(\TB_p)=0$ for $p\ge 0$ and $k \neq m$. In the case $k=m$, 
since $H^{m-1}_c(\TB_{p-1})=0$ and $H^{m+1}_c(\TB_p \setminus \TB_{p-1})=0$, 
the long exact sequence implies the result.
\endproof
\medskip

\textbf{Proof of Theorem \ref{thm:dimension}. }First we note 
that the space $\TB_m$ is a $2m$-dimensional 
topological manifold with boundary $\rd \TB_m=\TB_{m-1}$ 
(see Remark \ref{rem:corner}). 
By Lefschetz duality (\cite[Theorem 3.43 and Exercises 35]{Hat}), 
we have
\[
H_k(\TB_m,\TB_{m-1})=H_k(\TB_m,\rd \TB_m) \cong H^{2m-k}_c(\TB_m). 
\]
Hence $\dim_{\Q}H_k(\TB_m,\TB_{m-1})=\dim_{\Q}H^{2m-k}_c(\TB_m)=0$ 
for $k \neq m$ by Lemma \ref{lem:dim2}. 
In the case $k=m$, again by Lemma \ref{lem:dim2}, 
we have
\[
\dim_{\Q} H_m(\TB_m,\TB_{m-1})=\dim_{\Q} H^{m}_c(\TB_m)
=\sum_{p=0}^m \dim_{\Q}H_c^m(\TB_p \setminus \TB_{p-1}).
\]
By Poincar\'e duality and Lemma \ref{lem:dim1}, the dimension of 
$H_c^m(\TB_p \setminus \TB_{p-1})$ is given by 
\[
\dim_{\Q}H_c^m(\TB_p \setminus \TB_{p-1})=
\dim_{\Q} H_p(\TB_p \setminus \TB_{p-1})=
\binom{rn+m-p-1}{m-p} \binom{n+p-2}{p}.
\]
Finally, by Vandermonde's identity for multi-set coefficients, we have
\[
\sum_{p=0}^m \binom{rn+m-p-1}{m-p} \binom{n+p-2}{p}
=\binom{rn+n+m-2}{m}
\]
and this implies the result.
\endproof
\begin{cor}
We have an equality of $\Q$-vector spaces
\[
\L_{n,m}^{(r)}\otimes_{\R}\Q=\H_{n,m}^{(r)}\otimes_{\R}\Q.
\]
In particular, the standard multifork classes $\{\,[\mathbb{F}_{\k}] 
\,\vert\, \k \in K_{n,m}^{(r)}\,\}$ form a basis of 
the homological representation $\H_{n,m}^{(r)}$ over $\Q$. 
\end{cor}
\textbf{Proof. }It immediately follows from Proposition \ref{prop:standard} and 
Theorem \ref{thm:dimension}. 
\endproof
\begin{rem}
In this paper, we don's discuss the problem whether $\L_{n,m}^{(r)}$ coincides 
with $\H_{n,m}^{(r)}$ over $\R$. For details of this problem, 
we refer to \cite{Big3,PP}.
\end{rem}

\subsection{Framed Burau representations}
In this section, we introduce the framed Burau representation of 
the framed braid group, which contains the Burau representation of 
the braid group as the quotient representation. 
We also show that the reduced framed Burau representation coincides with 
the homological representation in the case $m=1$. 

Let $M_n^{(r)}$ be a free $\Z[q^{\pm}]$-module of rank $(n+rn)$ spanned 
by the basis
\[
\left\{\,a_i,b_i^{(s)}\,\middle\vert\,  
i=1,\dots,n \text{ and }s=1,\dots,r\,  \right\}. 
\]
Define the action of  $FB_n$ on  $M^{(r)}_n$ by
\begin{align*}
\sigma_i \cdot a_k&:=
\begin{cases}
\,(1-q)a_i+q a_{i+1}-(b_i^{(1)}+\cdots+b_i^{(r)}) &\text{ if }\:k=i \\
\,a_i   &\text{ if }\:k=i+1 \\
\,a_k &\text{ otherwise}
\end{cases}
\\
\sigma_i\cdot b_k^{(s)}&:=
\begin{cases}
\,q \,b_{i+1}^{(s)} &\text{ if }\:k=i \\
\,b_i^{(s)}   &\text{ if }\:k=i+1 \\
\,b_k^{(s)} &\text{ otherwise}
\end{cases}
\\
\tau_i \cdot a_k&:=
\begin{cases}
\,a_{i}-q\,b_i^{(r)} &\text{ if }\:k=i \\
\,a_k  &\text{ otherwise} \\
\end{cases}
\\
\tau_i \cdot b_k^{(s)}&:=
\begin{cases}
\,b_i^{(s-1)} &\text{ if }\:k=i \text{ and }\: s=2,\dots,r \\
\,q\,b_i^{(r)} &\text{ if }\:k=i \text{ and }\: s=1 \\
\,b_k^{(s)}   &\text{ otherwise}. 
\end{cases}
\end{align*}
\begin{prop}
The above action is well-defined. Hence it 
gives a representation of $FB_n$ on $M_{n}^{(r)}$.
\end{prop}
\textbf{ Proof. }By direct computation, we can 
easily check that the above action of 
$FB_n$ on $M_n^{(r)}$ satisfies the relations in Definition \ref{def:framed_braid}. 
\endproof
\medskip

We call $M_n^{(r)}$ the \textit{framed Burau representation} of 
the framed braid group $FB_n$. 
Let $N_n^{(r)} \subset M_n^{(r)}$ be the submodule spanned by 
$\{b_i^{(s)}\,\vert\,  i=1,\dots,n \text{ and }s=1,\dots,r\,  \}$. 
Then $N_n^{(r)}$ is a subrepresentation of $FB_n$. 
If we restrict the representation $N_n^{(r)}$ on the subgroup 
$B_n \subset FB_n$, then it 
becomes the direct sum of the standard representations 
of the braid group in \cite{TYM,Sys}. 
Thus our homological representation gives the homological interpretation 
of the standard representations of the braid groups. 
On the other hand, the quotient 
representation $M_n^{(r)} \slash N_n^{(r)}$ becomes 
the Burau representation of the braid group \cite{Bur} 
(see also \cite[Section 3.1]{KT}). 

Set $c_i:=a_{i+1}-q^{-1}a_i$ for $i=1,\dots,n-1$ and consider the 
submodule $L_n^{(r)}\subset M_n^{(r)}$ of rank $(n-1+rn)$ spanned by 
\[
\left\{c_1,\dots,c_{n-1},b_1^{(1)},\dots,b_1^{(r)},\dots,
b_n^{(1)},\dots,b_n^{(r)}\right\}. 
\]
Then it is easy to check that $L_n^{(r)}$ is closed under 
the action of $FB_n$. 
We call $L_n^{(r)}$ the \textit{reduced framed Burau representation} 
of the framed braid group $FB_n$. Note that $L_n^{(r)}$ also contains  
$N_n^{(r)}$ as a subrepresentation. Further, 
the quotient representation $L_n^{(r)} \slash N_n^{(r)}$ 
is isomorphic to the reduced Burau representation of 
the braid group (see \cite[Section 3.3]{KT}).

\begin{prop}
The reduced framed Burau representation $L_n^{(r)}$ is 
equivalent to the homological representation $\L_{n,1}^{(r)}$. 
\end{prop}
\textbf{ Proof. }Define $\k_i \in K_{n,1}^{(r)}$ and 
$\mathbf{l}_i^{(s)} \in K_{n,1}^{(r)}$ by
\begin{align*}
\k_i &:=(0,\dots,0,k_i=1,0,\dots,0)\\
\l_i^{(s)} &:=(0,\dots,0,l_i^{(s)}=1,0,\dots,0) .
\end{align*}
Then we can check that the action of $FB_n$ on fork classes
$[\mathbb{F}_{\k_i}]$ and $[\mathbb{F}_{\l_i^{(s)}}]$ 
coincides with the action on $c_i$ and $b_i^{(s)}$ respectively by 
using the fork rules in Figure \ref{Fig:Rules}. 
Hence the correspondence
\[
[\mathbb{F}_{\k_i}] \mapsto c_i,\quad 
[\mathbb{F}_{\l_i^{(s)}}] \mapsto b_i^{(s)}
\]
gives an isomorphism of representations $\L_{n,1}^{(r)} \iso L_n^{(r)}$. 
\endproof

\begin{rem}
\label{rem:Alexander}
Let $\psi_n \colon B_n \to \GL(n-1,\Z[q^{\pm 1}])$ be the reduced 
Burau representation of the braid group $B_n$. 
Then it is well-know that the quantity 
\[
\Delta_{\hat{\beta}}(q)=\frac{q-1}{q^n-1}\det (I_{n-1}-\psi_n(\beta))
\]
for $\beta \in B_n$ is the Alexander polynomial of 
the link $\hat{\beta}$ which is the plat closure of $\beta$
(see \cite[section 3.4]{KT}). 
Similarly, let
$\psi_n^{(r)} \colon FB_n \to \GL(n-1+rn,\Z[q^{\pm 1}])$ be 
the reduced framed Burau representation of the framed braid group 
$FB_n$ and consider the quantity
\[
\Delta_{\hat{\beta}}^{(r)}(q)=\frac{q-1}{q^n-1}
\det (I_{n-1+rn}-\psi_n^{(r)}(\beta))
\]
for $\beta \in FB_n$. Then it gives a polynomial invariant of 
the framed link $\hat{\beta}$. 
Details of this topic will be discussed in \cite{Ike}.
\end{rem}

\section{Monodromy representations}
\label{sec:monodromy}
In this section, we give a detailed study of the confluent KZ equations introduced 
in \cite{JNS} and construct representations of the framed braid groups 
as the monodromy representations of the confluent KZ equations. 
\subsection{Confluent Verma module bundles}
\label{sec:confluent_Verma}
Consider the Lie algebra $\g:=\mathfrak{sl}_2$ over $\C$ 
with the standard basis $\{E,H,F\}$ satisfying
\begin{align*}
[E,F]=H,\quad [H,E]=2E,\quad [H,F]=-2F.
\end{align*}
Define a \textit{truncated current Lie algebra} 
$\g^{(r)}$ by 
\[
\g^{(r)}:=\g[t] \slash t^{r+1}\g[t]
\]
where $\g[t]:=\g \otimes \C[t]$ is a Lie algebra with the bracket 
\[
[X \otimes t^m, Y \otimes t^n]:=[X,Y] \otimes t^{m+n}.
\]
Set $X_p:=X \otimes t^p$ for $X \in \g$ and $p=0,1,\dots,r$.
The Lie algebra $\g^{(r)}$ has the triangular decomposition
\[
\mathfrak{n}_+^{(r)}:=\bigoplus_{p=0}^r \C E_p,\quad 
\mathfrak{h}^{(r)}:=\bigoplus_{p=0}^r \C H_p,\quad
\mathfrak{n}_-^{(r)}:=\bigoplus_{p=0}^r \C F_p.
\]
Take complex numbers 
$\lam,\gamma^{(1)},\gamma^{(2)},\dots,\gamma^{(r)} \in \C$ 
with $\gamma^{(r)} \neq 0$ and set $\gamma:=(\gamma^{(1)},\dots,\gamma^{(r)})$. 
A \textit{highest weight vector $v_{\lam}(\gamma)$} is defined by
the conditions
\[
\mathfrak{n}_+^{(r)}  v_{\lam}(\gamma)=0,\quad 
H_0 v_{\lam}(\gamma)=\lam\, v_{\lam}(\gamma),\quad 
H_p  v_{\lam}(\gamma)=\gamma^{(p)}\, v_{\lam}(\gamma)
\]
for $p=1,\dots,r$. We call $\lam$  
a \textit{weight} and $\gamma^{(1)},\dots,\gamma^{(r)}$ 
\textit{movable weights}. 
Later, movable weights $\gamma^{(1)},\dots,\gamma^{(r)}$ are regarded as 
variables of the confluent KZ equation.    
\begin{defi}
A \textit{confluent Verma module} $M^{(r)}_{\lam}(\gamma)$ of a weight 
$\lam \in \C$, movable weights 
$\gamma=(\gamma^{(1)},\dots,\gamma^{(r)}) \in \C^{r-1}\times \C^*$ and 
\textit{P-rank} $r>0$ 
is defined to be the induced module
\[
M^{(r)}_{\lam}(\gamma):=
U\g^{(r)} \otimes_{U\mathfrak{b}^{(r)} }\C v_{\lam}(\gamma)
\]
where $\C v_{\lam}(\gamma)$ is the one dimensional representation of 
the Borel subalgebra 
$\mathfrak{b}^{(r)}:=\mathfrak{n}_+^{(r)} \oplus \mathfrak{h}^{(r)}$, 
and $U\g^{(r)}$ and $U\mathfrak{b}^{(r)}$ are 
universal enveloping algebras 
of $\g^{(r)}$ and $\mathfrak{b}^{(r)}$.  
\end{defi}
Later, we will see that the P-rank $r$ controls the Poincar\'e ranks 
of irregular singularities of the confluent KZ equation. 
It is known that $M^{(r)}_{\lam}(\gamma)$ is irreducible if $\gamma^{(r)} \neq 0$ 
(see \cite{Wil}). Let $\Z_+$ be non-negative integers. 
For $j:=(j^{(0)},j^{(1)},\dots,j^{(r)}) \in \Z_+^{r+1}$, write
\[
F^j v_{\lam}(\gamma):=F_0^{j^{(0)}}F_1^{j^{(1)}} \cdots 
F_r^{j^{(r)}}v_{\lam}(\gamma).
\]
Then the set $\{F^j v_{\lam}(\gamma)\,\vert\,j \in \Z_+^{r+1}\}$ forms a basis 
of $M^{(r)}_{\lam}(\gamma)$, and hence we have
\[
M^{(r)}_{\lam}(\gamma)\cong 
\bigoplus_{j \in \Z_+^{r+1}}\C F^j v_{\lam}(\gamma).
\]
As mentioned above, to treat movable weights $\gamma^{(1)},\dots,\gamma^{(r)}$  
as variables, we introduce the \textit{space of movable weights} by
\[
B^{(r)}:=\{\,(\gamma^{(1)},\dots,\gamma^{(r)})\in 
\C^r\,\vert\,\gamma^{(r)}\neq 0\,\} \cong \C^{r-1}\times \C^*. 
\]
Consider the family of confluent Verma modules over $B^{(r)}$
as follows. 
\begin{defi}
A \textit{confluent Verma module bundle} of a weight $\lam \in \C$ 
and P-rank $r>0$ is defined by
\[ 
E_{\lam}^{(r)}:=\bigcup_{\gamma \in B^{(r)}}M_{\lam}^{(r)}(\gamma)  
\to B^{(r)}
\]
where a fiber over $\gamma \in B^{(r)}$ is  
the confluent Verma module $M_{\lam}^{(r)}(\gamma) $.
\end{defi}

\subsection{Connections on confluent Verma module bundles}
\label{sec:connection}
As in the previous section, we use the coordinate
$(\gamma^{(1)},\dots,\gamma^{(r)}) \in B^{(r)} $ 
for the space of movable weights $B^{(r)}$.   
Introduce the holomorphic vector fields $D^{(0)},D^{(1)},\dots,D^{(r-1)}$ 
on $B^{(r)}$ by 
\[
D^{(s)}:=\sum_{p=1}^{r-s}\,p \,\gamma^{(s+p)} \frac{\rd}{\rd \gamma^{(p)}}. 
\]
Denote by $\O(B^{(r)})$ the space of holomorphic functions on $B^{(r)}$.
\begin{lem}
\label{lem:holomorphic}
The vector fields $D^{(0)},D^{(1)},\dots,D^{(r-1)}$
form a basis of the space of holomorphic vector fields on $B^{(r)}$ 
as an $\O(B^{(r)})$-module, 
i.e. any holomorphic vector filed $V$ on $B^{(r)}$ can be uniquely written as
\[
V=\sum_{s=0}^{r-1}f_s D^{(s)}
\]  
with $f_0,\dots,f_{r-1} \in \O(B^{(r)})$.
\end{lem}
\textbf{Proof. }First we note that 
$\{\rd \slash \rd \gamma^{(s)}\}_{s=1}^r$ is a basis of 
the space of holomorphic vector fields on $B^{(r)}$. 
Consider the base change between 
$\{D^{(s)}\}_{s=0}^{r-1}$ and 
$\{\rd \slash \rd \gamma^{(s)}\}_{s=1}^r$.  
The base change matrix
\[
\left(\frac{\rd}{\rd \gamma^{(1)}},\frac{\rd}{\rd \gamma^{(2)}},
,\dots \frac{\rd}{\rd \gamma^{(r)}}  \right)P=
(D^{(r-1)},D^{(r-2)},\dots,D^{(0)})
\]
is given by the upper triangular matrix
\[
P=
\begin{pmatrix}
\gamma^{(r)} & \gamma^{(r-1)} & \cdots & \gamma^{(1)} \\
                &2 \gamma^{(r)} & \cdots & 2\gamma^{(2)} \\
                &                    &      \ddots       & \vdots          \\
              &                   &             &   r \gamma^{(r)}
\end{pmatrix}.
\]
Since $\det P=r! (\gamma^{(r)})^r$ and it is nonzero on $B^{(r)}$, 
we can write $\rd \slash \rd \gamma^{(s)}$ as
\[
\frac{\rd}{\rd \gamma^{(s)}}=\frac{1}{(\gamma^{(r)})^r}
\sum_{p=1}^s f_{s,p}(\gamma) D^{(r-p)}
\]
with $f_{s,p}(\gamma) \in 
\C[\gamma^{(1)},\dots,\gamma^{(r)}]$. 
Hence the result follows. 
\endproof
\medskip

Let $\g^{(r)} \times B^{(r)} \to B^{(r)}$ be a trivial bundle 
and $\Sec(\g^{(r)} \times B^{(r)})$ be its holomorphic sections. 
\begin{defi}
For the vector fields $D^{(0)},\dots,D^{(r-1)}$, 
the connection
\[
\nabla_{D^{(s)}} \colon \Sec(\g^{(r)} \times B^{(r)})
\to \Sec(\g^{(r)} \times B^{(r)})
\]
is defined by
\[
\nabla_{D^{(s)}} \left( fX_p  \right)
:= (D^{(s)} f) X_p +f \,p X_{p+s}
\]
where $f \in \O(B^{(r)})$ 
and $X_p=X\otimes t^p$ for $X \in \g$. 
\end{defi}
\begin{lem}
\label{lem:hol_int}
The connection $\nabla$ is integrable. 
\end{lem}
\textbf{Proof. }By direct computation, we have
\[
[D^{(s)},D^{(t)}]=-(s-t)D^{(s+t)}
\]
and
\[
[\nabla_{D^{(s)}},\nabla_{D^{(t)}}] X_p=-(s-t) pX_{p+s+t}
=-(s-t) \nabla_{D^{(s+t)}}X_{p}.
\]
This implies the integrability
\[
[\nabla_{D^{(s)}},\nabla_{D^{(t)}}] -\nabla_{ [D^{(s)},D^{(t)}] }=
[\nabla_{D^{(s)}},\nabla_{D^{(t)}}] +(s-t)\nabla_{ D^{(s+t)}}=0.
\]
\endproof
\medskip

The computation in the proof of Lemma \ref{lem:holomorphic} 
implies that the connection $\nabla$ is 
a meromorphic connection on $\C^r$ with a pole 
along $\{\gamma^{(r)}=0\}$. It is a regular singularity if $r=1$ and  
an irregular singularity if $r \ge 2$. 
This connection can be extended to a connection on 
\[
U\g^{(r)} \times B^{(r)} \to B^{(r)}
\]
by defining
\[
\nabla_{D^{(s)}}  (U_1 U_2\cdots U_k):=
\sum_{l=1}^k\, \left( \nabla_{D^{(s)}} U_l \right)\,  U_1U_2\cdots 
\widehat{U_l} \cdots U_k 
\]
for $U_1,\dots,U_l \in \g^{(r)}$. It is well-defined since 
$\nabla$ preserves the relations in $U\g^{(r)}$, i.e.
\[
\nabla_{D^{(s)}}(U_1 U_2-U_2 U_1)=\nabla_{D^{(s)}} [U_1,U_2]
\]
for $U_1,U_2 \in \g^{(r)}$. 
As a result, we can also extend the connection on 
confluent Verma module bundles compatible with the action of $\g^{(r)}$. 
Let $E_{\lam}^{(r)} \to B^{(r)}$ be a confluent Verma module 
bundle. Denote by $v_{\lam}$ 
the global section of $E_{\lam}^{(r)}$ defined 
by gathering highest weight vectors 
$v_{\lam}(\gamma)$ from each fiber $M_{\lam}^{(r)}(\gamma)$. 
Since the rank of $E_{\lam}^{(r)}$ is infinite, we consider the space of 
holomorphic sections of $E_{\lam}^{(r)}$ as
\[
\Sec (E_{\lam}^{(r)}):=\left\{\, \sum_{\text{finite}} f_k\,
F^{j_k}v_{\lam}\,\middle\vert\, f_k \in \mathcal{O}(B^{(r)}) 
\text{ and }j_k \in \Z_+^{r+1}   \,\right\}.
\]
Then the connection $\nabla$ on 
$\g^{(r)} \times B^{(r)}$ induces a connection on $E_{\lam}^{(r)}$ 
as follows. 
\begin{defi} 
\label{def:connection}
For the vector fields $D^{(0)},\dots,D^{(r-1)}$, 
the connection
\[
\nabla_{D^{(s)}} \colon \Sec (E_{\lam}^{(r)}) \to \Sec (E_{\lam}^{(r)})
\]
is defined by
\[
\nabla_{D^{(s)}} (f F^j v_{\lam}):=(D^{(s)}f) F^j v_{\lam}
+f \,(\nabla_{D^{(s)}}F^j)v_{\lam}
\]
where $f \in \O(B^{(r)})$.  
\end{defi}
By Lemma \ref{lem:hol_int}, the connection $\nabla$ 
is integrable.

\subsection{Spaces of singular vectors}
Fix positive integers $r_1,\dots,r_n>0$. 
Let $M_{\lam_i}^{(r_i)}(\gamma_i)$ be 
confluent Verma modules of weights $\lam_i \in \C$ and movable weights 
$\gamma_i \in B^{(r_i)}$ for $i=1,\dots,n$. 
Set $\Lam:=(\lam_1,\dots,\lam_n)$, 
$\Gamma:=(\gamma_1,\dots,\gamma_n)$ and $R:=(r_1,\dots,r_n)$, and
consider the tensor product
\[
M^{(R)}_{\Lam}(\Gamma):=M_{\lam_1}^{(r_1)}(\gamma_1)\otimes \cdots \otimes 
M_{\lam_n}^{(r_n)}(\gamma_n).
\]
Then the direct sum of truncated current Lie algebras 
$\g^{(R)}:=\g^{(r_1)}\oplus \cdots \oplus \g^{(r_n)}$ 
naturally acts on  $M^{(R)}_{\Lam}(\Gamma)$.
The original Lie algebra $\g$ acts on $M^{(R)}_{\Lam}(\Gamma)$
through the diagonal embedding of $\g$ into $\g^{(R)}$.
In other words, an element $X \in \g$ acts 
on $M^{(R)}_{\Lam}(\Gamma)$ as 
\[
(X_0 \otimes 1 \otimes \cdots \otimes 1)+
(1 \otimes X_0 \otimes \cdots \otimes 1)+\cdots+
(1 \otimes 1 \otimes \cdots \otimes X_0).
\]
\begin{defi}
Set $|\Lam|:=\lam_1+\cdots+\lam_n$. 
For $m \in \Z_{\ge 0}$, define the  
\textit{space of vectors of weight} $(|\Lam|-2m)$ by 
\[
W^{(R)}_{\Gamma}[|\Lam|-2m]:=
\{\,w \in M^{(R)}_{\Lam}(\Gamma) \,\vert\,     
H w=(|\Lam|-2m )w\,\},
\] 
and the \textit{space of singular vectors of weight} 
$(|\Lam|-2m)$ by
\[
S^{(R)}_{\Gamma}[|\Lam|-2m]:=
\{\,w \in W_{\Gamma}^{(R)}[|\Lam|-2m]\,\vert\,     
E w=0\,\}.
\]
\end{defi}
Let $v_{\Lam}(\Gamma):=v_{\lam_1}(\gamma_1)\otimes \cdots \otimes 
v_{\lam_n}(\gamma_n)$ be the tensor product of highest weight vectors.
Set $
\Z_+^{(R+\one)}:=\Z_+^{r_1+1} \times 
\cdots \times \Z_+^{r_n+1}$. For $J=(j_1,\dots,j_n) \in \Z_+^{(R+\one)}$, write
\[
F^J v_{\Lam}(\Gamma):=F_0^{j_1^{(0)}}\cdots 
F_{r_1}^{j_1^{(r_1)}}v_{\lam_1}(\gamma_1)  \otimes 
\cdots \otimes F_0^{j_n^{(0)}}\cdots 
F_{r_n}^{j_n^{(r_n)}}v_{\lam_n}(\gamma_n) .
\]
Then the space $W^{(R)}_{\Gamma}[|\Lam|-2m]$
has the basis
\[
\left\{ \,F^J v_{\Lam}(\Gamma)\,\middle\vert\, J \in \Z_+^{(R+\one)} 
\text{ with } |J|=m\,  \right\}
\]
where $|J|:=\sum_{i=1}^n \sum_{s=0}^{r_i}j_i^{(s)}$. 
Thus the dimension of $W^{(R)}_{\Gamma}[|\Lam|-2m]$ is given by
\[
\dim_{\C} W^{(R)}_{\Gamma}[|\Lam|-2m] =\binom{|R|+n+m-1}{m}
\]
where $|R|=r_1+\cdots+r_n$. 
\begin{prop}
\label{prop:dimension}
The dimension of $S^{(R)}_{\Gamma}[|\Lam|-2m]$ is given by
\[
\dim_{\C} S^{(R)}_{\Gamma}[|\Lam|-2m] =\binom{|R|+n+m-2}{m}.
\]
\end{prop}
\textbf{Proof. }For simplicity, we write 
\[W^{(R)}:=W^{(R)}_{\Gamma}[|\Lam|-2m], 
\quad S^{(R)}:=S^{(R)}_{\Gamma}[|\Lam|-2m]
\]
and $x:=\gamma_1^{(r_1)}$. 
Set $R^{\prime}:=(r_1-1,r_2,\dots,r_n) $ and consider the subspace 
$W^{(R^{\prime})}\subset W^{(R)}$ spanned by
$\{\,F^J v_{\Lam}(\Gamma)\,\vert\, J \in \Z_+^{(R+\one)} 
\text{ with } |J|=m \text{ and }j_1^{(r_1)}=0\,  \}$. 
Then 
\[
\dim_{\C} W^{(R^{\prime})}=\binom{|R|+n+m-2}{m}.
\]
Define the operators $f,\,\rd_f$ and $e$ acting on $M^{(R)}_{\Lam}(\Gamma)$ 
by
\[
f:=F_{r_1}\otimes 1 \otimes \cdots \otimes 1,\quad
\rd_f:=\frac{\rd}{\rd F_{r_1}} \otimes 1 \otimes \cdots \otimes 1
\]
and $e:=E- x\rd_f$.
We show that the linear map $L\colon W^{(R)} \to W^{(R)}$ 
defined by
\[
L(v):=\exp(-fe \slash x)v= \sum_{k=0}^{\infty}\frac{(-1)^k (f e)^k v}{k! \,x^k  } 
\]
for $v \in W^{(R)}$ is well-defined and 
induces an isomorphism of subspaces $L \colon W^{(R^{\prime})} \iso S^{(R)}$. 
By direct computation, we can check that $[\,f \rd_f,e]=0$ and $[\,f \rd_f,fe]=fe$. 
Since 
\[
(f \rd_f) \,(fe)^k F^J v_{\Lam}(\Gamma)=(j_1^{(r_1)}+k)(fe)^k  F^J v_{\Lam}(\Gamma)
\]
and $j_1^{(r_1)} \le m$, we have $(f e)^k  F^J v_{\Lam}(\Gamma)=0$ for $k \gg 0$. 
Hence $L$ is well-defined. Next we show that $E\cdot L(v)=0$ if and only if 
$v \in W^{(R^{\prime})}$. Since $E=e+x \rd_f$, the condition $E\cdot L(v)=0$ 
implies that
\[
\rd_f \frac{(fe)^{k+1} v }{(k+1)!}= \frac{e(fe)^k v}{k!} 
\]
for all $k \ge 0$ and $\rd_f v=0$. We show it by induction. 
Note that $\rd_f v=0$ if 
and only if $v \in W^{(R^{\prime})}$. For $k=0$, 
\[
\rd_f(fev)=ev+f\rd_f (ev)=ev+ ef\rd_f v=ev
\]
by $[f \rd_f,e]=0$. For general $k$, 
\begin{align*}
\rd_f \frac{(fe)^{k+1} v}{(k+1)!}= \frac{e(fe)^{k} v}{(k+1)!}
+\frac{ f\rd_f \{e (fe)^k v\} }{(k+1)!}
= \frac{e(fe)^{k}v}{(k+1)!}+\frac{ e f\rd_f  \{(fe)^k v\} }{(k+1)\cdot k!}
\end{align*}
again by $[f \rd_f,e]=0$.
By induction hypothesis, 
\[
\rd_f \frac{(fe)^{k+1} v}{(k+1)!}=\frac{e(fe)^{k} v}{(k+1)!}+\frac{ef}{k+1}\cdot 
\frac{e (fe)^{k-1} v }{(k-1)!}
=\frac{e (fe)^k v}{k!}.
\]
Thus $L$ maps $W^{(R^{\prime})}$ into $S^{(R)}$. 
Clearly, $L$ is invertible.
\endproof

\subsection{Confluent KZ equations}
Let  $E^{(r_i)}_{\lam_i}\to B^{(r_i)} $
be confluent Verma module bundles of weights $\lam_i \in \C$ and 
P-rank $r_i>0$  
for $i=1,\dots,n$. 
Consider the external tensor product 
\[
E_{\Lam}^{(R)}:=E^{(r_1)}_{\lam_1}\boxtimes \cdots \boxtimes E^{(r_n)}_{\lam_n} 
\to B^{(r_1)} \times \cdots \times B^{(r_n)}
\]
where $\Lam:=(\lam_1,\dots,\lam_n)$ and $R:=(r_1,\dots,r_n)$. 
As the notation of coordinates on the base space 
$B^{(R)}:=B^{(r_1)}\times \cdots \times B^{(r_n)}$, we use
\[
(\gamma_1,\dots,\gamma_n) \in B^{(r_1)} \times \cdots \times B^{(r_n)},\quad
\gamma_i=(\gamma_i^{(1)},\dots,\gamma_i^{(r)}) \in B^{(r_i)} = \C^{r_i-1}\times \C^*.
\]
Define the elements $\Omega^{(p,q)}_{ij} \in U\g^{(r_i)} \otimes U\g^{(r_j)}$ by
\[
\Omega^{(p,q)}_{ij}:=E_p \otimes F_q+F_p \otimes E_q+\frac{1}{2}H_p \otimes H_q
\]
for $p=0,\dots,r_i$ and $q=0,\dots,r_j$. We also 
denote by $\Omega_{ij}^{(p,q)} \colon E_{\Lam}^{(R)} \to E_{\Lam}^{(R)}$
the induced bundle map through 
the action of $\Omega^{(p,q)}_{ij}$ on $i$-th and $j$-th components of $E_{\Lam}^{(R)}$. 
Let $
X_n:=\{\,(z_1,\dots,z_n) \in \C^n\,\vert\,\,z_i \neq z_j \text{ if }i \neq j  \,\} $
be the configuration space. 
We can extend the vector bundle
$E_{\Lam}^{(r)} \to B^{(R)}$ trivially on the direct product 
$ B^{(R)} \times X_n$ by  
the pullback of the projection $ B^{(R)}\times X_n \to B^{(R)}$.
By abuse of notation, we also write the resulting vector bundle as 
$E_{\Lam}^{(R)} \to B^{(R)} \times X_n$.
The following operators  
were introduced in \cite{JNS}.
\begin{defi}	
The operators
\[
\left\{\,G_i^{(s)} \in \End(E_{\Lam}^{(R)})\,\middle\vert\,
i=1,\dots,n \text{ and }s=-1,\dots,r_i-1\, \right\},
\]
called the \textit{(generalized) Gaudin Hamiltonians}, are defined by
\[
G_i^{(s)}:=\frac{1}{2}\sum_{p=0}^s \Omega_{i,i}^{(p,s-p)}+
\sum_{j \neq i}\sum_{\substack{p,q \ge 0 \\ p+q \le r_i+r_j-s-1}}
\binom{p+q}{p  }\frac{(-1)^p}{(z_i-z_j)^{p+q+1}}\,\Omega_{ij}^{(s+p+1,q)}. 
\] 
\end{defi}

Define the space of holomorphic sections 
of $E_{\Lam}^{(R)}$ by
\[
\Sec (E_{\lam}^{(R)}):=\left\{\, \sum_{\text{finite}} f_k\,
F^{J_k}v_{\Lam}  
\,\middle\vert\, f_k \in \mathcal{O}(B^{(R)}\times X_n) 
\text{ and }J_{k} \in \Z_+^{R+\one}   \,\right\}
\]
where $v_{\Lam}:=v_{\lam_1}\otimes \cdots \otimes v_{\lam_n}$. 
Then the action of Gaudin Hamiltonians on $\Sec (E_{\Lam}^{(R)})$
is well-defined.
As in Section \ref{sec:connection}, we introduce 
vector fields $D_i^{(s)}$ for $i=1,\dots,n$ and $s=0,\dots,r_i-1$ 
on $B^{(R)} \times X_n$ by
\[
D_i^{(s)}:=\sum_{p=1}^{r_i-s}\,p \,\gamma_i^{(s+p)} \frac{\rd}{\rd \gamma_i^{(p)}}. 
\]
Set $\rd_i:=\rd \slash \rd z_i$. 
Then
\[
\left\{\,\rd_i,\,D_i^{(s)}\,\middle\vert\,i=1,\dots,n 
\text{ and }s=0,\dots,r_i-1 \right\}
\]
forms a basis of the space of vector fields on $B^{(R)}\times X_n$. 
We extend the connection in Definition \ref{def:connection} 
as follows. 
\begin{defi}
The connection $\nabla$ on $E_{\Lam}^{(R)}$ is defined by
\begin{align*}
\nabla_{D_i^{(s)}} \left( f\, F^J v_{\Lam}\right) :=
\left(D_i^{(s)}f \right)F^J v_{\Lam}
+f \,F^{j_{1}}v_{\lam_1}\otimes 
\cdots \otimes \left(\nabla_{D_i^{(s)}} F^{j_i}v_{\lam_i}\right)  
\otimes \cdots \otimes F^{j_n}v_{\lam_n}
\end{align*}
and
\[
\nabla_{\rd_i}\left( f\,F^J v_{\Lam}\right)
:=\frac{\rd f}{\rd z_i}\,F^J v_{\Lam}
\]
where $f \in \O(B^{(R)}\times X_n)$.
\end{defi}
By Lemma \ref{lem:hol_int}, $\nabla$ is integrable. 
Now we introduce the confluent KZ equation as follows. 
\begin{defi}
Fix a nonzero complex number $\kappa \in \C^*$. 
The \textit{cofluent Knizhnik-Zamolodchikov (KZ) equation} for a 
holomorphic section $\Phi \in \Sec (E_{\Lam}^{(R)})$ 
is a system of differential equations
\begin{align*}
\nabla_{\rd_i} \Phi&=\frac{1}{\kappa}G_{i}^{(-1)} \Phi \\
\nabla_{D_i^{(s)}}\Phi&=\frac{1}{\kappa}
\left(G_{i}^{(s)} -\beta_i^{(s)}\right)\Phi
\end{align*}
for $i=1,\dots,n$ and $s=0,\dots,r_i-1$ where $\beta_i^{(s)}$ is 
the function given by
\[
\beta_i^{(s)}:=\frac{1}{4}\,\sum_{p=0}^s \gamma_i^{(p)}\gamma_i^{(s-p)}
+\frac{s+1}{2}\gamma_i^{(s)}.
\] 
\end{defi}
By the definition of Gaudin Hamiltonians, it has poles of order 
$r_i+r_j+1$ along the divisors $\{z_i=z_j\}$ for $i \neq j$. 
In addition, by the computation 
in Lemma \ref{lem:holomorphic},  
it has poles of order $r_i$  
along the divisors $\{\gamma_i^{(r_i)}=0\}$. Thus the confluent KZ equation has 
irregular singularities and their Poincar\'e ranks are determined by 
P-ranks $r_1,\dots,r_n$. 
In \cite{JNS}, they showed the integrability of the confluent KZ equation. 
\begin{prop}[\cite{JNS}, Proposition 4.1]
Define the confluent KZ connection $\nabla^{\KZ}$ by 
\begin{align*}
\nabla_{\rd_i}^{\KZ}&:=\nabla_{\rd_i}-\frac{1}{\kappa}
G_i^{(-1)}  \\
\nabla_{D_i^{(s)}}^{\KZ}&:=\nabla_{D_i^{(s)}}
-\frac{1}{\kappa}\left(G_{i}^{(s)} -\beta_i^{(s)}\right)
\end{align*}
for $i=1,\dots,n$ and $s=0,\dots,r_i-1$.
Then  $\nabla^{\KZ}$ is integrable.  
\end{prop}
Consider the restriction of the confluent KZ equation on 
finite rank subbundles of $E_{\Lam}^{(R)}$. We note the 
following property of the confluent KZ connection. 
\begin{lem}
\label{lem:commute}
The action of the Lie algebra $\g$ on $\Sec (E_{\Lam}^{(R)})$  
commutes with the action of the confluent KZ connection $\nabla^{\KZ}$ 
on $\Sec (E_{\Lam}^{(R)})$.
\end{lem}
\textbf{Proof. }It is easy to check that $[\,\Omega_{ij}^{(p,q)}, X]=0$ for 
$X \in \g$. In addition, $[\,\nabla,X]=0$ by the definition of $\nabla$.  
Hence $[\nabla^{\KZ}, X]=0$.
\endproof
\medskip

Let $S^{(R)}[|\Lam|-2m]$ be the finite rank subbundle of $E_{\Lam}^{(R)}$ 
whose fiber over $(\Gamma,z) \in B^{(R)}\times X_n$ is 
the space of singular vectors $S^{(R)}_{\Gamma}[|\Lam|-2m]$. 
Lemma \ref{lem:commute} implies that the restriction of 
the confluent KZ equation on the subbundle $S^{(R)}[|\Lam|-2m]$ 
is well-defined. 

\subsection{Monodromy representations}
In this section, we construct the monodromy 
representations of the framed braid groups 
from the confluent KZ connections. 
We first note that the space of movable weights $B^{(r)}=\C^{r-1}\times \C^*$ 
is homotopic to the circle $S^1$. 
Hence, the space $B^{(R)}=B^{(r_1)}\times \cdots \times B^{(r_n)}$ 
is homotopic to the $n$-dimensional torus $T^n$ 
and the direct product $B^{(R)} \times X_n$ is homotopic to $T^n \times X_n$. 
Thus the fundamental group  $\pi_1(B^{(R)} \times X_n,*)$ is 
isomorphic to the pure 
framed braid group $FP_n$ (see Section \ref{sec:fundamental}). 

Recall from the previous section that 
the confluent KZ connection is an integrable connection on 
the vector bundle $E_{\Lam}^{(R)} \to B^{(R)} \times X_n$. 
Consider the restriction of the confluent KZ connection on 
the finite rank subbundle $S^{(R)}[|\Lam|-2m] \subset E_{\Lam}^{(R)}$. 
Then we obtain the monodromy representation 
of the pure framed braid group
\[
\theta_{\Lam,\kappa}^{(R)} \colon 
FP_n  \to \Aut_{\C} S^{(R)}_{\Gamma}[|\Lam|-2m]
\] 
where $S_{\Gamma}^{(R)}(\Gamma)[|\Lam|-2m]$ is the fiber over 
the base point $*=(\Gamma,z) \in B^{(R)} \times X_n$. 
The dimension of the representation 
$\theta_{\Lam,\kappa}^{(R)}$ is given 
by Proposition \ref{prop:dimension}.

Consider the case $R=(r,\dots,r)=(r^n)$ and $\Lam=(\lam,\dots,\lam)=(\lam^n)$. 
In this case, the action of the symmetric group $\mathfrak{S}_n$ 
on $B^{(r^n)} \times X_n$ is well-defined since $B^{(r^n)}=(B^{(r)})^n$.
In addition,  $\lam_1=\cdots =\lam_n$ implies that 
the action of $\mathfrak{S}_n$ lifts on $E_{\lam^n}^{(r^n)}$. 
Hence we obtain the quotient vector bundle
\[
E_{\lam^n}^{(r^n)} \slash \mathfrak{S}_n 
\to (B^{(r^n)} \times X_n )\slash \mathfrak{S}_n.
\]
Further, the confluent KZ connection $\nabla^{\KZ}$ 
is $\mathfrak{S}_n$-invariant in the case $r_1=\cdots=r_n$. 
Thus the connection $\nabla^{\KZ}$ descends to an 
integrable connection on the quotient vector bundle 
$E_{\lam^n}^{(r^n)} \slash \mathfrak{S}_n $. 
By Proposition \ref{prop:fundamental}, the fundamental group 
$\pi_1((B^{(r^n)} \times X_n )\slash \mathfrak{S}_n,*)$ is 
isomorphic to the framed braid group $FB_n$. 
Therefore we obtain the following monodromy representation. 
\begin{prop}
For positive integers $r,m>0$ and complex parameters 
$(\lam,\kappa) \in \C \times \C^*$, 
there is a representation of the framed braid group 
\[
\theta_{\lam,\kappa}^{(r)} \colon FB_n \to \Aut_{\C} S^{(r^n)}_{\Gamma}[n\lam-2m]
\]
constructed as the monodromy representation of the 
confluent KZ connection on the vector bundle
\[
S^{(r^n)}[n \lam-2m] \slash \mathfrak{S}_n 
\to (B^{(r^n)} \times X_n) \slash \mathfrak{S}_n.
\] 
\end{prop}
The dimension of this representation is given by Proposition \ref{prop:dimension}.

\end{document}